\documentclass[12pt]{article} 
\usepackage[numbers]{natbib}
\usepackage{array,epsfig,fancyheadings,rotating}

\usepackage{amsmath}
\usepackage{amssymb}
\usepackage{amsfonts}
\usepackage{amsthm}

\usepackage{graphicx}
\usepackage{bbm}
\usepackage{float}
\usepackage{sectsty}
\newcommand*{\Scale}[2][4]{\scalebox{#1}{$#2$}}

\sectionfont{\large}
\subsectionfont{\normalsize}

\RequirePackage[OT1]{fontenc}
\RequirePackage{amsthm,amsmath}

\numberwithin{equation}{section}
\newtheorem{thm}{Theorem}

\newtheorem{cor}{Corollary}

\theoremstyle{definition}

\newtheorem{re}{Remark}
\newtheorem{df}{Definition}

\begin{document}

\markboth{\hfill{\footnotesize\rm Tadeusz Inglot, Teresa Ledwina and Bogdan \'Cmiel} \hfill}
{\hfill {\footnotesize\rm INTERMEDIATE EFFICIENCY} \hfill}

\renewcommand{\cite}[1]{%
  \citeauthor{#1}~(\citeyear{#1})}

\renewcommand{\thefootnote}{}
$\ $\par


\fontsize{12}{14pt plus.8pt minus .6pt}\selectfont \vspace{0.8pc}
\centerline{\large\bf INTERMEDIATE EFFICIENCY} 
\centerline{\large\bf IN NONPARAMETRIC TESTING PROBLEMS}
\vspace{2pt} \centerline{\large\bf WITH AN APPLICATION TO}
\vspace{2pt} \centerline{\large\bf  SOME WEIGHTED STATISTICS}
\vspace{.4cm} \centerline{Tadeusz Inglot$^*$, Teresa Ledwina$^\diamond$ and Bogdan \'Cmiel$^{\diamond,\dagger}$} \vspace{.4cm} \centerline{\it
Wroc{\l}aw University of Science and Technology $^*$ } 
\centerline{\it  Polish Academy of Sciences $^\diamond$ and  AGH University of Sciences and Technology $^\dagger$} \vspace{.55cm} \fontsize{9}{11.5pt plus.8pt minus
.6pt}\selectfont

\begin{center}
{\it Dedicated to Wilbert Kallenberg, with friendship and esteem}
\end{center}


\begin{quotation}
\noindent {\bf Abstract.}
 The basic motivation and primary goal of this paper is a qualitative evaluation of the performance of a new weighted statistic for a nonparametric test for stochastic dominance based on two samples, which was introduced in \cite{r24}. For this purpose, we elaborate a useful variant of Kallenberg's notion of intermediate  efficiency. This variant is general enough to be applicable to other nonparametric problems.
We provide a formal definition of the proposed variant of intermediate efficiency, describe the technical tools used in its calculation, and provide proofs of related asymptotic results. Next, we apply this approach to calculating the intermediate efficiency of the new test with respect to the classical one-sided Kolmogorov-Smirnov test, which is a recognized standard  for this problem. 
It turns out that for a very large class of convergent alternatives the new test is more efficient than the classical one. 
We also report the results of an extensive simulation study on the powers of the  tests considered, which shows that the new variant of intermediate efficiency reflects the exact behavior of the power well. \par

\vspace{9pt}
\noindent{\bf Mathematics Subject Classification}. 62G10, 62G20, 60E15.\\

\noindent {\it Key words and phrases.} Anderson-Darling weight, 
asymptotic relative efficiency,
Kallenberg efficiency, 
Kolmogorov-Smirnov test,
local alternatives,
moderate deviations,
rank empirical process,
stochastic dominance,
stochastic order,
two-sample problem.

\end{quotation}\par

\fontsize{12}{14pt plus.8pt minus .6pt}\selectfont

\section{Introduction}\label{S1}
\fancyhf{}
\rhead{\thepage}
\lhead{\centerline{INTERMEDIATE EFFICIENCY}}
\rfoot{} 

In order to assess the performance of a test, a multitude of concepts of efficiency have been proposed. 
See \cite{r32} and \cite{r37} for an overview of earlier definitions of asymptotic relative efficiencies.
More recently, some efficiency measures have been defined in terms of probabilities of 
 large and moderate deviations  of type I and type II errors of tests; cf. \cite{r39}, \citeauthor{r40} (\citeyear{r40,r41}), \cite{r44} for some ideas, related history and implementations. Obviously, each approach has its own rationale, some limitations, inherent complexity, and was developed with its own specific assumptions. 

In this paper, we concentrate on the concept of asymptotic relative efficiency (ARE), and restrict our attention to Kallenberg's notion of intermediate efficiency.  We develop a variant of this notion which is useful in nonparametric problems and enables widely applicable, tractable, analytic comparisons between two tests. The results from earlier applications of Kallenberg's concept of efficiency are very encouraging. \cite{r22} gave several implementations of this notion to some popular tests for selected  parametric and semiparametric models. Further developments include: adaptive test statistics (\cite{r13}, \citeauthor{r18} (\citeyear{r18,r19})), some classical goodness-of-fit statistics for both continuous and discrete data  (\citeauthor{r19} (\citeyear{r19,r20}), and \cite{r47}), some classical statistics for testing a simple parametric hypothesis for a model of one-sample censorship (\cite{r45}), testing for no-effect in certain regression models (\cite{r20}, \cite{r46}).

In Section \ref{S3} we apply our variant of this approach to a recent solution of the one-sided nonparametric test for stochastic dominance
introduced in \cite{r24}. Section 5.4.3 of that paper contains an extensive discussion motivating such an investigation. In fact, here we consider a generalization of that solution and denote it by ${\cal T}_N$.
The test  statistic ${\cal T}_N$ is asymptotically equivalent to the maximum of a weighted rank empirical process over a grid in (0,1), where the end points  depend on the sample size and the classical Anderson-Darling weight is used. This last statistic is denoted by ${\cal W}_N$. We compare the new solution with the classical (unweighted) two-sample Kolmogorov-Smirnov test ${\cal V}_N$, which is commonly applied to this problem. In order to obtain a formula for the efficiency, several technical results have to be proved.
Bounds for the asymptotic power of ${\cal T}_N$, under sequences of alternatives, and  moderate deviations under the appropriate null distributions, are the main technical results of this paper.
We prove that, for this one-sided nonparametric test, carefully matching the weight and the range of the maximum in ${\cal W}_N$ is highly profitable and results in a test which dominates the classical unweighted solution, regardless of how the alternative deviates from the null hypothesis. In particular, our Theorem 6 shows that the intermediate efficiency of the new test with respect to the Kolmogorov-Smirnov test is greater or equal than 1 for a very large class of sequences of alternatives.  In simple terms, we compare two consistent tests which, in principle, can detect any fixed alternative. The only question is how many observations are needed for them to attain a given power. The intermediate efficiency of a weighted test relative to an unweighted one is the number which indicates by which factor one has to increase the sample size when using the unweighted solution to have approximately the same power as the weighted procedure has. Our simulations illustrate that the value of this efficiency measure appropriately explains the empirical behavior of these tests' power in this sense.

We close the introduction with a discussion of Kallenberg's notion of ARE and some related problems. The efficiency measure was introduced by matching the basic features of the approaches of Pitman and Bahadur to ARE. This involves alternatives converging to the null model (slower than in Pitman's approach) and significance levels tending to 0 (slower than in Bahadur's theory).  
Obviously, according to this approach, both the significance level and the alternative depend on the sample size $n$. The essence of this setting is that the significance level goes to 0 as $n$ increases, while the asymptotic power, under the underlying sequence of local alternatives, should be non-degenerate. Hence, the sum of the type I and type II errors is in (0,1).
Such requirements call for a careful and delicate balancing of the rates at which the significance levels tend to 0 and the local alternatives approach the null model. 
Although Kallenberg's concept of efficiency is slightly more complicated than the classical notions of Pitman and Bahadur, it is more widely applicable. 

The advantages and limitations of Pitman's approach are well known; see \cite{r32} and \cite{r38}  for insightful comments. In particular, in most non-elementary applications, when the underlying test statistic is not asymptotically normal, then Pitman's efficiency may not exist or may depend on the significance level and a given power. 
Bahadur's efficiency is much more widely applicable than Pitman's ARE. However, it requires the existence of non-degenerate large deviations for the test statistics under consideration and this turns out to be  impossible to guarantee for many statistics used at present, including several weighted ones; cf. \cite{r5}, \cite{r42}, and \cite{r43}.
In comparison with Bahadur's ARE, the concept of intermediate efficiency requires similar, but less demanding conditions, to be applicable. 
\cite{r41}, p. 622, admits: ``Kallenberg's efficiency represents the analogue of the Bahadur efficiency in a moderate deviation zone''. 

The paper of \cite{r22} is mathematically elegant. Its main part concerns the one-sample case. In addition, it puts emphasis on having results which hold for all possible sequences of local alternatives.
These rather stringent conditions can result in it being impossible to implement, even in relatively simple situations.
Moreover, one important question regarding ensuring the non-degeneracy of the asymptotic power was skipped in all the examples given in that paper, 
while the $k$-sample case was not described in detail. \cite{r13}  provided  further analysis  of Kallenberg's efficiency in the one-sample case, mainly in the context  of studying adaptive Neyman tests. 
Our Remark \ref{re5} in Section \ref{S2} carefully discusses the differences between Kallenberg's original approach and our contribution.

In view of the above, in the present paper we propose a simple as possible  variant of the notion of intermediate efficiency and define tools which help to calculate it. In contrast to the original concept, we do not require that results hold for all local alternative sequences.  Our focus is on relaxing the requirements on the test statistics applied by as much as possible. Moreover, we embed the one- and two-sample cases into a joint scheme and study them simultaneously. 
Our contribution is self-contained: all of the necessary tools are carefully described, sufficient conditions for the non-degeneracy of the limiting power are given, 
and thorough proofs of all the related results are provided. The details are presented in Section \ref{S2}, as well as Appendices A and B.

The organization of the paper is as follows: Section \ref{S2} describes our setup, the pathwise variant of intermediate efficiency, 
 and discusses previous work in the light of our proposals. After this preparatory material,
Section \ref{S3} presents two constructions of tests (a new and a classical one) for detecting stochastic ordering, introduces a natural class of sequences of alternatives for this problem, gives the theoretical results needed to compare such tests via intermediate ARE, and gives an explicit formula for the intermediate efficiency of the new solution with respect to the classical Kolmogorov-Smirnov test. We conclude by reporting some results of an extensive simulation study, which aims to show the usefulness of this variant of ARE when evaluating powers for finite samples. 
The proofs of all the technical results are given in Appendices A - G.  Appendix B also provides some technical lemmas, which are useful in checking the assumptions of Theorem 1, the main theorem.  \\

\section{Intermediate efficiency. Pathwise variant. Basic facts and comments}\label{S2}

\subsection{Notation and definitions}\label{S2.1}
Consider two independent samples $X_1,...,X_m$ and $Y_1,...,Y_n$ defined on the same measurable space and coming from the probability distributions $P$ and $Q$, respectively. The situation $P=Q$ corresponds to the one-sample case. Generalization to the $k$-sample case and independence testing are obvious. 

Assume that $m=m(N),\; n=n(N)$, with $N=m+n$. In the one-sample case, one can assume that the sample size is $N$, $m=m(N)=\lfloor N/2 \rfloor$ and $n=N-m$. 
Suppose that $\eta_N=m(N)/N$ satisfies $\eta_N \to \eta \in (0,1),$ as $N \to \infty$. Denote the set of all product measures $P\times Q$ under consideration by  $\mathbb{P}$. Let $\mathbb{P}_0$ be a subset of $\mathbb{P}$. We want to test whether 
$$
\mathbb{H}_0 : P\times Q \in \mathbb{P}_0\;\;\;\mbox{or}\;\;\;\;\mathbb{H}_1 : P\times Q \in \mathbb{P}_1 = \mathbb{P} \setminus \mathbb{P}_0
$$
is true. 

Suppose that we have two upper-tailed tests defined by sequences of real valued test statistics ${\cal U}^{(I)}_N$ and ${\cal U}^{(II)}_N$
and critical values corresponding to a significance level $\alpha \in (0,1)$ and sample size $N$, denoted by $u^{(I)}_{\alpha N}$ and $u^{(II)}_{\alpha N}$, respectively. To be specific, let
$$
u^{(I)}_{\alpha N} = \inf\Bigl\{w : \sup_{P\times Q \in \mathbb{P}_0} P^{m(N)}\times Q^{n(N)} \bigl({\cal U}^{(I)}_N > w\bigr) \leq \alpha\Bigr\},
$$
where $P^{m(N)}$ and $Q^{n(N)}$ are $m(N)$ and $n(N)$ fold products of $P$ and $Q$, respectively. We define $u^{(II)}_{\alpha N}$ in an analogous way. Throughout this article, $\{s_N\}$ denotes an infinite sequence of elements $s_N$, $N\geq 1$. 

We wish to evaluate the efficacy of the test based on  ${\cal U}^{(II)}=\{{\cal U}^{(II)}_N\}$ by comparing its sensitivity relative to the sensitivity of the test based on ${\cal U}^{(I)}=\{{\cal U}^{(I)}_N\}$. Hence, the test based on ${\cal U}^{(I)}$ is used as a benchmark. 

We shall consider significance levels which tend to 0 as the sample size grows. Namely, the set $\mathbb{L}$ of sequences of all admissible levels in the intermediate setting is defined
by
\begin{equation}\label{2.1}
\mathbb{L}= \bigl\{\{\alpha_N\} : \alpha_N \to 0,\;\;N^{-1} \log \alpha_N \to 0\bigr\}.
\end{equation}

Note that (\ref{2.1}) excludes significance levels $\alpha_N$ which tend to 0 exponentially fast. This is a characteristic of Bahadur efficiency.

Given $P_0 \times Q_0 \in \mathbb{P}_0$, $P_1 \times Q_1 \in \mathbb{P}_1$, and $\vartheta \in (0,1)$, define
\begin{equation}\label{2.2}
P_{\vartheta}=(1-\vartheta)P_0 + \vartheta P_1,\;\;\;\;\;\;Q_{\vartheta}=(1-\vartheta)Q_0 + \vartheta Q_1,
\end{equation}
and assume that 
\begin{equation}\label{2.3}
P_{\vartheta} \times Q_{\vartheta} \in \mathbb{P}_1\;\;\;\mbox{for every}\;\;\vartheta \in (0,1).
\end{equation}
Hereafter, the measures $P_0 \times Q_0 \in \mathbb{P}_0$ and  $P_1 \times Q_1 \in \mathbb{P}_1$ defining $P_{\vartheta} \times Q_{\vartheta}$ and satisfying (\ref{2.3}) are fixed.

Now, consider a particular sequence $\{\theta_N\}$, $\theta_N \in (0,1)$, where $\theta_N \to 0,\;$as $N \to \infty$, and the corresponding sequence $\{P_{\theta_N}\times Q_{\theta_N}\}$. Given $N$ and the corresponding $\theta_N$, for an arbitrary natural number $S$ representing the total  size of an auxiliary sample and corresponding $m(S)$ and $n(S)$, define 
$\Pi_{\theta_{N}}^S\hspace{-0.1cm}=P_{\theta_N}^{m(S)}\times Q_{\theta_{N}}^{n(S)}.$
Thus $\Pi_{\theta_N}^N = P_{ \theta_N}^{m(N)} \times Q_{ \theta_N}^{n(N)}$. Finally, set 
$$\Scale[1.6]{\pi}=\left\{\Pi_{\theta_N}^N\right\}.$$
Suppose that there exists $\{\alpha_N\}=\{\alpha_N (\Scale[1.6]{\pi})\} \in \mathbb{L}$ which satisfies
\begin{equation}\label{2.4}
0 < \liminf_{N \to \infty} \Pi_{ \theta_N}^N \bigl({\cal U}_N^{(II)} > u^{(II)}_{ \alpha_N N}\bigr) \leq
\limsup _{N \to \infty} \Pi_{ \theta_N}^N \bigl({\cal U}_N^{(II)} > u^{(II)}_{ \alpha_N N}\bigr) < 1.
\end{equation}
Let 
\begin{equation}\label{2.5}
\mathbb{L}^{*} =  \mathbb{L}^* (\Scale[1.6]{\pi})=\Bigl\{\{\alpha_N\}=\{\alpha_N(\Scale[1.6]{\pi})\} \in \mathbb{L} : (\ref{2.4}) \mbox{ holds}\Bigr\}.
\end{equation}
In consequence, given $\{\theta_N\}$, for all $\{\alpha_N\} \in  \mathbb{L}^*$, the corresponding test based on ${\cal U}_N^{(II)}$ has non-degenerate asymptotic power under $\{P_{\theta_N} \times Q_{\theta_N}\}$. In the sequel, we assume that $\mathbb{L}^*$ is nonempty.

The definition of the intermediate efficiency of ${\cal U}^{(II)}$ with respect to ${\cal U}^{(I)}$, which we give below, refers to this particular sequence of alternatives $\{P_{\theta_N} \times Q_{\theta_N}\}$, together with the corresponding $\Scale[1.6]{\pi}$, and sequences of significance levels $\{\alpha_N\} \in \mathbb{L}^*.$ 
For every $N\geq1$, let
\begin{equation}\label{2.6}
M_{{\cal U}^{(II)} {\cal U}^{(I)}}\bigl(N,\Scale[1.6]{\pi}\bigr)=
\end{equation}
\vspace{-0.5cm}
\begin{equation}\nonumber
\inf\Bigl\{M \geq 1 : \Pi_{ \theta_N}^N \bigl({\cal U}_N^{(II)} > u^{(II)}_{\alpha_N N}\bigr) \leq
\Pi_{\theta_N}^{M+k} \bigl({\cal U}_{M+k}^{(I)} > u^{(I)}_{\alpha_N {M+k}}\bigr)\;\;\mbox{for all}\;k \geq 0\Bigr\}.
\end{equation}

\vspace{2mm}
\begin{df}
If 
\begin{equation}\label{2.7}
e_{{\cal U}^{(II)} {\cal U}^{(I)}}=
\lim_{N \to \infty}\frac{M_{{\cal U}^{(II)} {\cal U}^{(I)}}\bigl(N,\Scale[1.6]{\pi}\bigr) }{N} \in [0,\infty]
\end{equation}
exists and does not depend on the choice of $\{\alpha_N\} \in \mathbb{L}^*$, we say that the asymptotic intermediate efficiency of ${\cal U}^{(II)}$ with respect to ${\cal U}^{(I)}$, under the sequence of alternatives $\{P_{\theta_N} \times Q_{\theta_N}\}$, exists and equals $e_{{\cal U}^{(II)} {\cal U}^{(I)}}$.
\end{df}

Obviously, the asymptotic behavior of $M_{{\cal U}^{(II)} {\cal U}^{(I)}}\bigl(N,\Scale[1.6]{\pi}\bigr)$, and hence $e_{{\cal U}^{(II)} {\cal U}^{(I)}}$, depends on $\eta$. Typically, the value of $\eta\in (0,1)$ is fixed and therefore, to simplify the notation, this parameter is omitted in (\ref{2.7}). However, in Section 4 we analyze numerically the influence of $\eta$ on the efficiency of tests for the two-sample case. Therefore, in Theorem \ref{thm6}, which presents an analytic formula for this measure of efficiency, and in Section 4 we clearly indicate the dependence of efficiency on $\eta$.

\begin{re}\label{re1}
In contrast to the original definition of \cite{r22} and its extension by \cite{r13}, where a counterpart of (\ref{2.7}) was required for some families of sequences of alternatives, the above definition is restricted to a particular sequence. Hence, it can be considered to be a kind of pathwise variant of the previous approaches. The path $\Scale[1.6]{\pi}$ is uniquely determined via $P_0 \times Q_0$, $P_1 \times Q_1$, $\{\theta_N\}$, $\{m(N)\}$ and $\{n(N)\}$. Such a pathwise approach extends the range of possible applications of the notion and allows us to avoid many non-trivial technicalities. 
In particular, we avoid the introduction of so-called renumerable families, which are  key objects in \cite{r13}. Note also that,
as a rule, (\ref{2.7}) holds for many sequences simultaneously, but the above definition treats each of them separately. For an illustration of this, see Theorem \ref{thm6} and the comment following it.
\end{re}

\begin{re}\label{re2}
Since $\eta_N \to \eta \in (0,1)$, we can succinctly rephrase the interpretation of $e_{{\cal U}^{(II)} {\cal U}^{(I)}}$ as follows: the test corresponding to ${\cal U}^{(I)}$ and based on the sample sizes \\ $(\lfloor m e_{{\cal U}^{(II)} {\cal U}^{(I)}}\rfloor ,\lfloor n e_{{\cal U}^{(II)} {\cal U}^{(I)}}\rfloor)$ has approximately the same power, under 
$\{P_{\theta_N}\times Q_{\theta_N}\}$, as the power of the test corresponding to $\;{\cal U}^{(II)}$  and based on the sample sizes $(m,n)$. 

Another useful interpretation of the intermediate efficiency is the value of the shift of a non-parametric alternative  necessary for the two tests under consideration to have the same local power. See \cite{r19} for some simple illustration and \cite{r14} for further development.

\end{re} 
\begin{re}\label{re3} 
To prove that the limit in (2.7) exists
and to obtain an explicit formula for it, we need to introduce some regularity assumptions for both test statistics. Sequences of statistics satisfying assumptions of this kind are called Kallenberg sequences by \cite{r45}.
Similarly to \cite{r13} and in contrast to \cite{r22}, we impose stronger requirements on ${\cal U}^{(I)}$ than on ${\cal U}^{(II)}$. On one hand, the benchmark, ${\cal U}^{(I)}$, can be always chosen in a convenient way. On the other hand, any relaxation of the requirements on ${\cal U}^{(II)}$ extends the scope of possible applications of this approach. Similar to the Bahadur efficiency, the intermediate efficiency of ${\cal U}^{(II)}$ with respect to ${\cal U}^{(I)}$ is calculated as the ratio between two slopes. These slopes are determined by an index for moderate deviations under the null hypothesis and a scaling factor which results from a kind of weak law of large numbers (WLLN) under the sequence of alternatives. It is worth emphasizing that we only assume a knowledge of moderate deviations of ${\cal U}^{(II)}$ in some restricted range. 
\end{re}

\subsection{Regularity assumptions on ${\cal U}^{(I)}$}\label{S2.2} 
\noindent
(I.1) There exists  a positive number $c_{{\cal U}^{(I)}}$ such that for every sequence $w_N > 0$ satisfying $w_N \to 0$ and $Nw_N^2 \to \infty$, the following  holds: 
$$
-\lim_{N \to \infty} \frac{1}{Nw_N^2} \log \sup_{P \times Q \in \mathbb{P}_0} P^{m(N)} \times Q^{n(N)} \bigl({\cal U}^{(I)}_N > w_N\sqrt N\bigr) =c_{{\cal U}^{(I)}}.
$$

\noindent
(I.2) There exists a function $b_{{\cal U}^{(I)}}(P_{\vartheta}\times Q_{\vartheta}),\;\vartheta\in (0,1)$, and a number $\rho \in [1,2]$, such that for every sequence $\{\vartheta_N\}$ of positive numbers where $\vartheta_N\to 0$ and
$N\vartheta_N^{\rho}\to\infty$, then for the corresponding $\Pi_{\vartheta_N}^N=P_{\vartheta_N}^{m(N)}\times Q_{\vartheta_N}^{n(N)}$, and every $\epsilon >0$
\begin{equation}\label{2.8}
\lim_{N\to\infty} \Pi_{\vartheta_N}^N\left(\left|\frac{{\cal U}_N^{(I)}}{\sqrt{mn/N}b_{{\cal U}^{(I)}}(P_{\vartheta_N}\times Q_{\vartheta_N})}-1\right|\geq \epsilon\right)=0. 
\end{equation}

We call $c_{{\cal U}^{(I)}}$ the index of moderate deviations of ${\cal U}^{(I)}$, while the quantity 
$c_{{\cal U}^{(I)}}\bigl[b_{{\cal U}^{(I)}}(\Pi_{\vartheta_N}^N)\bigr]^2$, where $b_{{\cal U}^{(I)}}(\Pi_{\vartheta_N}^N)=\sqrt{\frac{mn}{N}}b_{{\cal U}^{(I)}}(P_{\vartheta_N} \times Q_{\vartheta_N})$,
is called the intermediate slope of ${\cal U}^{(I)}$ under $\Pi_{\vartheta_N}^N$.

\subsection{Regularity assumptions on ${\cal U}^{(II)}$}\label{S2.3}
\noindent
(II.1) There exist  sequences $\{\gamma_N\}$ and $\{\lambda_N\}$,
$1 \leq \gamma_N < \lambda_N \leq N$,  such that $\gamma_N/\lambda_N \to 0$ and for every positive sequence $w_N$, where $Nw_N^2/\lambda_N \to 0$ and $Nw_N^2/\gamma_N \to \infty$, the following holds 
\begin{equation}\label{2.9}
-\lim_{N \to \infty} \frac{1}{Nw_N^2} \log \sup_{P \times Q \in \mathbb{P}_0} P^{m(N)} \times Q^{n(N)} \bigl({\cal U}^{(II)} > w_N\sqrt N\bigr) =c_{{\cal U}^{(II)}}
\end{equation}
for a positive number $c_{{\cal U}^{(II)}}$.

\noindent
(II.2) For any sequence $\{\theta_N\}$  such that $\;\theta_N \to 0$ and
$N\theta_N^{\rho} \to \infty$, where $\rho$ is defined in (I.2), there exists a positive sequence 
$\{b_{{\cal U}^{(II)}}(\Pi_{\theta_N}^N)\}$  such that 
\begin{equation}\label{2.10}
\lim_{N \to \infty}  \Pi_{ \theta_N}^N \Bigl(\Bigl|\frac{{\cal U}^{(II)}_N}{b_{{\cal U}^{(II)}}(\Pi_{\theta_N}^N)} - 1\Bigr| \geq \epsilon\Bigr)=0
\;\;\;\;\mbox{ for every}\;\;\epsilon > 0. 
\end{equation}

As above, the quantity 
$c_{{\cal U}^{(II)}}\bigl[b_{{\cal U}^{(II)}}(\Pi_{\theta_N}^N)\bigr]^2$ is called the intermediate slope of ${\cal U}^{(II)}$ under $\Pi_{\theta_N}^N$.

\noindent
\begin{re}\label{re4} 
In (I.2) and (II.2) we imposed the assumption $N \theta_N^{\rho} \to \infty,$ for some $\;\rho \in [1,2]$. In several previously considered cases,  it was simply assumed that $\rho =2$. As seen from the proof of Theorem \ref{thm1}, given in Appendix A, this assumption is closely related to the behavior of the probability of a type one error in the  test under consideration,
which in turn depends explicitly on the intermediate slope. Obviously, the most demanding conditions are provided by the Neyman-Pearson test. The intermediate slope of the Neyman-Pearson test for two-sample problems was studied in \cite{r9}; cf. Theorems 3.3 and 3.4. In this case, $\rho=2$ is appropriate. In other tests, the situation might be different. 
\end{re}

\subsection{Computation of the intermediate efficiency}\label{S2.4} 
We start with a counterpart to Lemmas 2.1 and 5.1 in \cite{r22}. This result gives conditions under which the asymptotic ratio of the slopes of ${\cal U}^{(I)}$ and ${\cal U}^{(II)}$ coincides with the asymptotic ratio of the sample sizes in (\ref{2.7}). Next, we compare our result to the lemmas mentioned above and comment on how the assumptions of our result can be checked.

\noindent
\begin{thm}\label{thm1}
Assume that for a given ${\cal U}^{(I)}$ the conditions (I.1) and (I.2) are true for some $\rho \in [1,2]$. Suppose that, for a particular sequence $\{\theta_N\}$, ${\cal U}^{(II)}$ satisfies (II.2). Moreover, (II.1) holds for some given $\{\gamma_N\}$ and $\{\lambda_N\}$. Suppose that there exists $\{\alpha_N\}$ from $\mathbb{L}^*$ such that
\begin{equation}\label{2.11}
\lim_{N \to \infty} \frac{\log \alpha_N}{\lambda_N} = 0\;\;\;\;\mbox{and}\;\;\;\lim_{N \to \infty} \frac{\log  \alpha_N}{\gamma_N} = -\infty.
\end{equation}
Finally, assume that the following limit exists
\begin{equation}\label{2.12}
\lim_{N\to\infty}\frac{c_{{\cal U}^{(II)}}[b_{{\cal U}^{(II)}}(\Pi_{\theta_N}^N)]^2}{c_{{\cal U}^{(I)}}(mn/N)[b_{{\cal U}^{(I)}}(P_{\theta_N}\times Q_{\theta_N})]^2}={\bf e} \in [0,\infty].
\end{equation}
Then the intermediate efficiency (\ref{2.7}) of $\;{\cal U}^{(II)}$ with respect to $\;{\cal U}^{(I)}$, under a particular sequence of alternatives 
$\{P_{\theta_N} \times Q_{\theta_N}\}$, exists and $e_{{\cal U}^{(II)} {\cal U}^{(I)}}={\bf e}$.
\end{thm}
Let us start by discussing the differences between our requirements and those in Kallenberg's Lemmas 2.1 and 5.1.

\begin{re}\label{re5}  

The first essential difference between these two papers consists in the fact that \cite{r22} requires that for both statistics ${\cal U}_N^{(I)}$ and ${\cal U}_N^{(II)}$ the WLLN under the alternatives holds for all allowable sequences $\{\theta_N\}$ and all corresponding distributions, provided that his condition (1.3) holds. This can be an inappropriate formulation and his paper provides examples where some further restrictions are imposed, and, in fact, some paths introduced.

The second essential difference between both papers lies in the fact that \cite{r22} assumes the same type of moderate deviations for both ${\cal U}_N^{(I)}$ and ${\cal U}_N^{(II)}$, and does not impose any restriction on the rate of convergence of $\{w_N\}$ (in our notation) from above. In contrast, we introduce such a requirement in the case of ${\cal U}_N^{(II)}$, and this is expressed by the assumption that $N w_N^2/\gamma_N \to \infty$. We impose such an assumption because it naturally arises when studying the  statistic ${\cal T}_N$; cf. Theorem \ref{thm3}. Such an assumption is simply indispensable in many cases, as, for example,  for some weighted goodness of fit statistics where moderate deviations exist and are only non-zero for a very restricted type of sequences of $\{w_N\}$'s. To balance this useful restriction on moderate deviations for ${\cal U}_N^{(II)}$, we assume that ${\cal U}_N^{(I)}$ has non-zero moderate deviations for the whole range of the $w_N$'s. This is not a restrictive assumption, because any convenient benchmark ${\cal U}_N^{(I)}$ can be chosen.

The third difference is that we do not require that $b_{{\cal U}^{(II)}}(\Pi_{\theta_N}^N)$ has a special structure. This allows us to compare statistics with a different rate of convergence than the one corresponding to ${\cal U}_N^{(I)}$.

Additionally, our paper clarifies the intermediate approach for the two-sample case. Note that, putting aside the limitations discussed above, the formulation of (iii) in Kallenberg's Lemma 5.1 contains some further non-explicit restrictions. Moreover, observe that the essential assumption regarding the non-degeneracy of the asymptotic power of ${\cal U}_N^{(II)}$, in our notation, is missing in the formulations of Kallenberg's lemmas, though it is clearly stated on p. 171 of \cite{r22}. This may be misleading, as it is needed in the proof. Verification of this assumption is a non-trivial problem, which is ignored in the analysis of the examples in \cite{r22}.  In our Theorem \ref{thm1}, this assumption is implicit in the condition $\{\alpha_N\} \in \mathbb{L}^*$, and later we provide some convenient tools for checking this assumption; cf. Remark \ref{re6} . 
\end{re}

Now,  we give some brief comments on verifying the assumptions of our Theorem 1. 

\begin{re}\label{re6}
The regularity assumptions (I.1) and (II.1) hold for many classical statistics. In the present paper, the two-sample Kolmogorov-Smirnov statistic serves as a good example. Often, the strong approximations method proves to be very useful in obtaining such results.

For some statistics, conditions like (I.2) and (II.2), together with the form of $\alpha_N$ and non-degeneracy of asymptotic powers, have already been justified on the basis of the limiting distribution of the underlying test statistic under given sequences of alternatives; cf. \citeauthor{r18} (\citeyear{r18,r19}) for some examples. However,  such a limiting distribution is often hard to derive. 
Therefore, following the idea applied for the first time in \cite{r20}, we propose to use some asymptotic bounds for the distributions of test statistics under the considered sequence of alternatives. These bounds have the condition of non-degenerate asymptotic powers explicitly built into them, which is very useful in checking (I.2) and (II.2)  in some relatively complex cases. For details, see Lemmas 1 and 2 in Appendix B, where we also give an example of a sequence $\{\alpha_N\}\in \mathbb{L}^*$ and a range of non-degenerate asymptotic powers for which Theorem 1 applies. See also Section 3.4 for an illustration of how this approach works.
\end{re}

\section{The intermediate efficiency of some Kolmogorov-Smirnov-type tests for stochastic ordering}\label{S3}
In this section, step by step we shall develop the tools necessary to apply the results of Section \ref{S2} to some selected statistics for testing for the existence of stochastic dominance. In successive subsections we present this problem and the test statistics under consideration, introduce appropriate sequences of alternatives, collect some theoretical results on moderate deviations under the null hypothesis and asymptotic behavior under the sequences of alternatives. We conclude with Theorem \ref{thm6}, which states the existence and the form of the intermediate efficiency, and Corollary \ref{cor1}, which exemplifies an implementation of the general result.

\subsection{The testing problem and local sequences of nonparametric alternatives}\label{S3.1}

We consider two independent samples $X_1,...,X_m$ and $Y_1,...,Y_n$ which correspond to the continuous distribution functions $F$ and $G$, respectively. As in Section \ref{S2.1}, we assume  that $m=m(N)$,  $n=n(N),\;m+n=N,$ and both sample sizes tend to infinity as $N$ increases. Moreover, we set $\eta_N = m(N)/N$ and suppose that 
\begin{equation}\label{3.1}
\eta=\lim_{N \to \infty} \eta_N \;\;\;\mbox{ exists and}\;\;\;  \eta \in (0,1).
\end{equation} 
Throughout Section 3, all limits are taken under the assumption that $m$ and $n$ grow in such a way that (\ref{3.1}) holds.

The null hypothesis,  $\mathbb{H}_0$, asserts that the $X$'s are stochastically smaller than the $Y$'s, i.e. 
$$
\mathbb{H}_0 : F(z) \geq G(z)\;\;\;\mbox{for each}\;\;\;z \in \mathbb{R},
$$
while the alternative, $\mathbb{H}_1$, is unrestricted and is of the form
$$
\mathbb{H}_1 : F(z) < G(z)\;\;\;\mbox{for some}\;\;\;z \in \mathbb{R}.
$$
Note that both classical tests for $\mathbb{H}_0$ and the new one, which we shall consider in the following subsections, are distribution free for any continuous $F=G$.
Moreover, the following property holds: 
$$
Pr({\cal S}_N > w | F \geq G) \leq Pr({\cal S}_N > w | F = G),\;\;\;\mbox{for all}\;\;w \in \mathbb{R},
$$
where ${\cal S}_N$ is any of the considered test statistics. For details, see \cite{r24}.  Therefore, to illustrate the approach based on intermediate efficiency with minimal technicality, we restrict our attention to
$F=G$, which satisfies $\mathbb{H}_0$.

To define a family of nonparametric paths, we proceed as in Remark 2.4 in \cite{r31}; see also \cite{r48} and \cite{r9}. Namely, we take an arbitrary $(F_1,G_1)$ satisfying $\mathbb{H}_1$ and select $F=G=F_0=J_1$ from $\mathbb{H}_0$, where
$
J_1(z)=\eta F_1(z) + (1-\eta)G_1(z),\;\;\;z \in \mathbb{R}.
$
The probability measure corresponding to this $F_0$ is denoted by $P_0$. $P_0^N$ denotes its $N$-fold product. 

With the above choice of $F_0$, for a given sequence of real numbers $\vartheta_N \in (0,1)$ such that $\vartheta_N \to 0$ as $N \to \infty$, we introduce a contamination model $(F_{1N},G_{1N})$  based on
\begin{equation}\label{3.2}
(F_{1N},G_{1N}) = (1-\vartheta_N)(F_0,F_0) + \vartheta_N(F_1,G_1).
\end{equation}
Thus $F_{1N}$ and $G_{1N}$ are absolutely continuous with respect to $J_1$ and possess densities $1-\vartheta_N (1-\eta) d(G_1-F_1)/dJ_1$ and 
$1+(\vartheta_N \eta)d(G_1-F_1)/dJ_1$, respectively. 
We denote by $P_{\vartheta_N}$ and $Q_{\vartheta_N}$ the probability measures corresponding to $F_{1N}$ and $G_{1N}$, defined in (\ref{3.2}). Formula (\ref{3.2}) defines our local sequence of nonparametric alternatives.

\subsection{One-sided two-sample test statistics}\label{S3.2}

The empirical distribution functions for the two samples are
$
\hat F_m(z) \! = \! {m}^{-1}\! \sum_{i=1}^m \! {\bf{1}}(X_i \leq z)$ and $\hat G_n(z) = {n}^{-1}\sum_{i=1}^n {\bf{1}}(Y_i \leq z)
$, respectively,
where  ${\bf 1}(\mathbb{E})$ is the indicator function of the set $\mathbb{E}$. Additionally, let $\hat J_N (z)$ be the empirical distribution function for the pooled sample, i.e.
\begin{equation}\label{3.3}
\hat J_N (z) = \eta_N \hat F_m (z) + (1-\eta_N)\hat G_n (z),\;\;\;z \in \mathbb{R},
\end{equation}
and denote the inverse of $\hat J_N (z)$ by $\hat J_N ^{-1}(t),\;t \in (0,1)$.
The one-sided Kolmogorov-Smirnov test rejects $\mathbb{H}_0$ when
\begin{equation}\label{3.4}
{\cal V}_N = \sqrt{\frac{mn}{N}} \sup_{z \in \mathbb{R}} \Bigl\{\hat G_n(z) - \hat F_m(z)\Bigr\}=
\sqrt{\frac{mn}{N}} \max_{1 \leq j \leq N} \Bigl\{\hat G_n - \hat F_m\Bigr\} \circ \hat J_N ^{-1}(\frac{j}{N})
\end{equation}
exceeds the appropriate critical value. Note that ${\cal V}_N $ is a rank statistic.

For testing $\mathbb{H}_0$ against $\mathbb{H}_1 $, \cite{r24} introduced, among other statistics, a test statistic based on the minimum of an appropriate set of linear rank statistics. Since $\mathbb{H}_0$ is one-sided, linear rank statistics with non-increasing score generating functions are appropriate; cf. \cite{r49}. In \cite{r24} step functions, related to  projections of Haar functions, were used to define a set of useful rank statistics. Here, we  slightly generalize this construction by allowing a more flexible set of step functions. 

To be specific, let $R_i,\;i=1,...,m,$ denote the rank of $X_i$ in the pooled sample $X_1,...,X_m,Y_1,...,Y_n$. Analogously,  $R_i,\;i=m+1,...,N,$ stands for the rank of $Y_i$ in the pooled sample. Set
\begin{equation}\label{3.5}
{\ell}_j(t) = {\ell}_{jN}(t)=- \sqrt{\frac{1-\pi_{jN}}{\pi_{jN}}}\, {\bf 1}(0\leq t <\pi_{jN}) +
\sqrt{\frac{\pi_{jN}}{1- \pi_{jN}}}\, {\bf 1}(\pi_{jN} \leq t \leq 1),
\end{equation}
where $0 < \pi_{1N} < \pi_{2N} < ... < \pi_{\Delta(N)N} < 1$, while $\Delta(N)$ is a non-decreasing sequence of natural numbers, such that $1 < \Delta(N) \leq N$.
Consider the corresponding linear rank statistics given by 
\begin{equation}\label{3.6}
{\cal L}_j ={\cal L}_{jN}= \sum_{i=1}^{N} c_{Ni}\,{\ell}_j\Bigl(\frac{R_i-0.5}{N}\Bigr),
\end{equation}
where
\begin{equation}\label{3.7}
c_{Ni} = \sqrt{\frac{mn}{N}}
\left\{\begin{array}{lrl}
- m^{-1} & \mbox{if} &  1\leq i \leq m, \\
\quad n^{-1} & \mbox{if} &  m <   i \leq N.
\end{array} \right.
\end{equation}
From the above, after elementary calculations, we obtain
$$
{\cal L}_j=\sqrt{\frac{mn}{N}} \frac{1}{\sqrt{\pi_{jN}(1-\pi_{jN})}}\times
$$ 
$$
\Bigl[\frac{1}{m} \sum_{i=1}^m {\bf 1}_{[0,\pi_{jN})}\Bigl(\frac{R_i-0.5}{N}\Bigr)
-\frac{1}{n} \sum_{i=m+1}^N {\bf 1}_{[0,\pi_{jN})}\Bigl(\frac{R_i-0.5}{N}\Bigr)\Bigr]
$$
$$
=\sqrt{\frac{mn}{N}} \frac{1}{\sqrt{\pi_{jN}(1-\pi_{jN})}} \Bigl[\frac{N}{mn} \sum_{i=1}^m {\bf 1}_{[0,\pi_{jN})}\Bigl(\frac{R_i-0.5}{N}\Bigr)-
\frac{\lceil N \pi_{jN} - 0.5 \rceil}{n}\Bigr],
$$
where $\lceil \bullet \rceil$ equals the number $\bullet$ when it is an integer, otherwise $\lceil \bullet \rceil = \lfloor \bullet \rfloor + 1$.
The above statistic, ${\cal L}_j$, contrasts the average values of the rescaled ranks $(R_i-0.5)/N,\;j=1,...,N,$ from the two samples which fall into 
the interval $[0,\pi_{jN})$. 
Equivalently, given the vector of ranks $(R_1,...,R_N)$, the value of ${\cal L}_j$ is a linear function of the values of the rescaled ranks $(R_i-0.5)/N,\;i=1,...,m,\;$from the first sample which fall into the interval $[0,\pi_{jN})$. Therefore, it is intuitive that small values of ${\cal L}_j$ indicate $\mathbb{H}_1$.

Finally, based on the union-intersection principle, set $M_{\Delta(N)} = \min_{1 \leq j \leq \Delta(N)} {\cal L}_j$. For a thorough discussion of the construction and properties of $M_{\Delta(N)}$ in the case where the $\pi_{jN}$'s are related to a dyadic partition of (0,1), we refer the reader to \cite{r24}. See also \cite{r25} for some useful properties of the ${\cal L}_j$'s.

Since upper-tailed critical regions have a long tradition  in efficiency calculations, for a given $\Delta(N)$,  we shall consider
\begin{equation}\label{3.8}
{\cal T}_N = \max_{1 \leq j \leq \Delta(N)} \bigl\{-{\cal L}_j\bigr\}
\end{equation}
and the related test which rejects $H_0$ for large values of ${\cal T}_N$. 

The statistic ${\cal T}_{N}$  differs by some asymptotically negligible quantity from the following weighted Kolmogorov-Smirnov-type statistic:
\begin{equation}\label{3.9}
{\cal W}_N = \sqrt{\frac{mn}{N}} \max_{1 \leq j \leq \Delta(N)} \frac{( \hat G_n - \hat F_m) \circ \hat J_N^{-1}(\pi_{jN})}{\sqrt{\pi_{jN} (1- \pi_{jN})}}.
\end{equation}
More precisely, the following holds for sufficiently large $N$:
\begin{equation}\label{3.10}
|{\cal T}_N - {\cal W}_N| \leq C_{\eta} /\sqrt {N \min\{\pi_{1N},1-\pi_{\Delta(N)N}\}},
\end{equation}
where $C_{\eta}$ is a positive number which depends only on $\eta$. For a proof, see Appendix D.

Our theoretical results are stated under the following basic assumptions: $\Delta(N)=o(N)$ and $1/[\Delta(N) +1] \leq \pi_{1N} < \pi_{\Delta(N)N} \leq 1 - 1/[\Delta(N) +1]$. Therefore, the right hand side of (\ref{3.10}) is $o(1)$ and any result that is proven  for ${\cal T}_N$ is automatically valid for ${\cal W}_N$ and vice versa. This gives us the flexibility to use the most appropriate and interpretable techniques for particular proofs. For convenience, all of the results are only stated for ${\cal T}_N$.

When applying the results of Section \ref{S2} to the above statistics, we set ${\cal U}_N^{(I)}={\cal V}_N$ and ${\cal U}_N^{(II)}={\cal T}_N$.

\subsection{Moderate deviations of ${\cal V}_N$ and ${\cal T}_N$ under $P^N_0$}\label{S3.3} 
Under $P_0^N$, defined in Section \ref{S3.1}, one can expect that the tails of ${\cal V}_N$ behave similarly to the tails of the classical Kolmogorov-Smirnov statistic for uniformity. Indeed, the ideas developed in \cite{r16} can also be applied to ${\cal V}_N$ and we obtain the following result:

\noindent
\begin{thm}\label{thm2} 
For any real sequence $\{w_N\}$ such that $w_N \to 0$ and $Nw_N^2 \to \infty$, the following holds:
\begin{equation}\label{3.11}
-\lim_{N \to \infty} \frac{1}{Nw_N^2} \log P_0^N\bigl({\cal V}_N \geq w_N \sqrt N \bigr) = c_{\cal V} = 2.
\end{equation}
\end{thm}
The constant $c_{\cal V} = 2$, appearing in (\ref{3.11}), is  the index of moderate deviations of ${\cal V}$.

To obtain moderate deviations result for  ${\cal T}_N$, we proceed as follows. From (\ref{3.8}), ${\cal T}_N$ is the maximum of $\Delta(N)$ linear rank statistics with non-continuous score functions, the ${\ell}_j'$s. These score functions are simple and can be approximated sufficiently well by piecewise linear functions. Some results proved in \cite{r15} can be applied to the corresponding rank statistics. 

\noindent 
\begin{thm}\label{thm3} 
Assume the following: 
(i) $\Delta(N) \to \infty,\;\;\Delta(N)=o(N),\;\; {as}$ $\;N \to \infty,$
(ii) $\pi_{1N} \geq 1/[\Delta(N)+1],\;\; \pi_{\Delta(N)N} \leq 1-1/[\Delta(N)+1],$ 
(iii) $\{w_N\}$ {\it is a sequence of positive numbers such that} $w_N \to 0,\;\;Nw_N^2/\log N \to \infty,\;$ {\it and,  for some} $\upsilon \in (0,1),\;\; w_N^{
1-\upsilon}\Delta(N) \to 0$ {\it as} $N \to \infty$. 
Then
\begin{equation}\label{3.12}
-\lim_{N \to \infty} \frac{1}{Nw_N^2} \log P_0^N\bigl({\cal T}_N \geq w_N \sqrt N \bigr) = c_{\cal T} = \frac{1}{2}.
\end{equation}
\end{thm}

From (\ref{3.12}), the index of moderate deviations of ${\cal T}$ is $c_{\cal T} =1/2.$

\subsection{The asymptotic behavior of ${\cal V}_N$ and ${\cal T}_N$ under  alternatives satisfying (\ref{3.2})}\label{S3.4}
Set 
$
\bar A(t) =\bar A(t;\eta)=(G_1-F_1)\circ J_1^{-1}(t),\;$ for $t\in (0,1).
$
Recall that, in this case, $P_{\vartheta_N}\sim F_{1N},\;Q_{\vartheta_N}\sim G_{1N},\;\Pi_{\vartheta_N}^N=P_{\vartheta_N}^{m(N)}\times Q_{\vartheta_N}^{n(N)},$ and ${\cal U}^{(I)}={\cal V}_N$. We shall show that (I.2) holds for ${\cal V}_N$ with 
$$
b_{{\cal U}^{(I)}}(P_{\vartheta} \times Q_{\vartheta})=\vartheta \sup_{0<t<1}\bar A(t;\eta)=\vartheta \sup_{z \in \mathbb{R}}[G_1(z)-F_1(z)].
$$
To simplify the notation, we also introduce
\begin{equation}\label{3.13}
b_{\cal V} (\Pi_{\vartheta_N}^N) = \sqrt{\frac{mn}{N}} b_{{\cal U}^{(I)}}(P_{\vartheta_N} \times Q_{\vartheta_N})
\end{equation}
and note that since $(F_1,G_1)$ belongs to $\mathbb{H}_1$, then $b_{\cal V} (\Pi_{\vartheta_N}^N) > 0$. The following result, along with Lemma 2 stated  in Appendix B, allows us to check whether ${\cal V}_N$ satisfies (I.2).

\noindent
\begin{thm}\label{thm4} 
Assume that $\vartheta_N \in (0,1),\; \vartheta_N \to 0$ and $N \vartheta_N^2 \to \infty$. Then
\begin{equation}\label{3.14}
\limsup_{N \to \infty} \Pi_{\vartheta_N}^N \bigl({\cal V}_N - b_{\cal V} (\Pi_{\vartheta_N}^N) \leq w\bigr) \leq V_2(w),\;\;\;w \in \mathbb{R},
\end{equation}
and
\begin{equation}\label{3.15}
\liminf_{N \to \infty} \Pi_{\vartheta_N}^N \bigl({\cal V}_N - b_{\cal V} (\Pi_{\vartheta_N}^N) \leq w\bigr) \geq V_1(w),\;\;\;w \in \mathbb{R}_+,
\end{equation}
where
$V_1(w) = Pr \bigl(\sup_{0 < t <1} B(t) \leq w\bigr),\;B$ is a Brownian bridge,
$V_2(w)=\Phi(\frac{w}{\sqrt{J_1(z_0)[1-J_1(z_0)]}}),$ and $\Phi$ denotes the $N(0,1)$ distribution function,
while $z_0= \inf\bigl\{z \in \mathbb{R} : G_1(z)-F_1(z) = \sup_{w \in \mathbb{R}}[G_1(w)-F_1(w)\bigr]\}$.
\end{thm}

Now take a particular sequence $\{\theta_N\}$ such that $\theta_N \in (0,1), \theta_N \to 0$ and $\;N\theta_N^2 \to \infty$ and the corresponding $\Pi_{\theta_N}^N$. Set 
$$
b_{\cal T} (\Pi_{\theta_N}^N) = \theta_N \sqrt{\frac{mn}{N}} \max_{1 \leq j \leq \Delta (N)} \frac{\bar A(\pi_{jN})}{\sqrt{\pi_{jN}(1-\pi_{jN})}}.
$$
Our next result gives conditions on $\Delta(N)$ and further restrictions on $\theta_N$ which guarantee that ${\cal U}_N^{(II)}={\cal T}_N$ satisfies (i) and (ii) of Lemma 1, given in Appendix B.

\noindent 
\begin{thm}\label{thm5} 
Suppose that the following conditions are satisfied: 
(i) $\pi_{1N} \geq 1/[\Delta(N)+1],\;\; \pi_{\Delta(N)N} \leq 1-1/[\Delta(N)+1]\;$ and 
$\max_{1 \leq j \leq \Delta(N)}  \bigl\{\pi_{jN} - \pi_{j-1 N}\bigr\} \to 0,$  when $\;\;N \to \infty$, 
(ii) $\theta_N \in (0,1),\; \theta_N \to 0$,  $\;N \theta_N^2 /\log^2 N  \to \infty$,  and  $\;\theta_N \Delta(N) \to 0, \;$ as $N \to \infty,$
(iii) $\eta_N \to \eta,$ $\eta\in (0,1)\;\;\;\mbox{and}\;\;\theta_N(\eta_N-\eta)\sqrt N = O(1),\;\;\;\mbox{when}\;\;N \to \infty$. 
Then
\begin{equation}\label{3.16}
\limsup_{N \to \infty} \Pi_{\theta_N}^N \bigl({\cal T}_N -  b_{\cal T} (\Pi_{\theta_N}^N) \leq w\bigr) \leq T_2(w),\;\;\;w \in \mathbb{R},
\end{equation}
and
\begin{equation}\label{3.17}
\liminf_{N \to \infty} \Pi_{\theta_N}^N \bigl({\cal T}_N -  b_{\cal T} (\Pi_{\theta_N}^N) \leq w\bigr) \geq T_1(w),\;\;\;w \in \mathbb{R}_+,
\end{equation}
where
$T_2(w)=\Phi(w)$, while $T_1(w)=Pr(\sup_{t \in [\delta,1-\delta]}|B(t)|/\sqrt{t(1-t)} \leq w)$ with $B$ being a Brownian bridge and $\delta = \delta(\bar A) \in (0,1)$  being defined by formula (G.3) in Appendix G. 
\end{thm} 

\subsection{The main result on the efficiency of ${\cal T}_N$ with respect to ${\cal V}_N$}\label{S3.5}  
The theoretical results presented above allow us to apply Theorem \ref{thm1} and to
formulate the following result on the intermediate efficiency of a class of tests based on ${\cal T}_N$. Recall that
$\overline{A}(t;\eta)=(G_1-F_1)\circ J_1^{-1}(t),\;\;J_1(z)=\eta F_1(z)+(1-\eta)G_1(z)$ and let
$A^*(t;\eta)= \overline{A}(t;\eta)/{\sqrt {t(1-t)}}$.

\begin{thm}\label{thm6}
Assume that conditions (i) and (iii) of Theorem \ref{thm5} hold and sharpen condition (ii) to (ii)' $\theta_N \in (0,1),\;\;N\theta_N^2 /\log^2N \to \infty,$ {\it and there exists} $\nu \in (0,1)$ {\it such that} $\theta_N^{1-\nu}\Delta(N) \to 0$ {\it when} $N \to \infty.$
Then, given $\eta \in (0,1)$, the   
intermediate efficiency $e_{{\cal T}{\cal V}}$ of ${\cal T}$ with respect to ${\cal V}$, under $\{P_{\theta_N} \times Q_{\theta_N}\}$, exists and is equal to 
\begin{equation}\label{3.18}
e_{{\cal T}{\cal V}}(\eta) = \frac{1}{4} \Bigl[\frac{\sup_{0<t<1}A^*(t;\eta)}{\sup_{0 < t < 1} \bar A(t;\eta)}\Bigr]^2
\end{equation}
$$=\frac{1}{4}\Bigl[\frac{\sup_{z \in \mathbb{R}}[(G_1(z)-F_1(z)]/\sqrt{J_1(z)[1-J_1(z)]}}{\sup_{z \in \mathbb{R}}\bigl[G_1(z)-F_1(z)\bigr]}\Bigr]^2.
$$
Moreover, for any sequence of alternatives it follows that 
$e_{{\cal T}{\cal V}}(\eta) \geq 1$ with equality holding if and only if $A^*(t;\eta)$ 
attains its maximum at $t=1/2$.
\end{thm}

\noindent
\begin{re}\label{re7}
It is worth emphasizing that in Theorem \ref{thm6} there are no restrictions on the deviation of the alternative $(F_1,G_1)$ from $(F_0,G_0)$. The assumptions concern only the rates of convergence of $\theta_N$ to $0$ and $\Delta(N)$ to $\infty$.
\end{re}

Theorem \ref{thm6} explains qualitatively the outcomes of the simulations in \citeauthor{r24} (\citeyear{r24,r26}). In Section 4 we demonstrate that the value of $e_{{\cal T}{\cal V}}(\eta)$ also gives precise quantitative information on the  relation between the empirical powers of sequences of the statistics ${\cal T}$ and ${\cal V}$ in the sense described in Remark \ref{re2}.

\noindent
\begin{re}\label{re8}
Theorem \ref{thm6} shows how the requirements on the rate of convergence of $\theta_N \to 0$ and $\Delta(N) \to \infty$ should be balanced in order to obtain the efficiency $e_{{\cal T}{\cal V}}(\eta)$. Corollary \ref{cor1}, stated below, gives a simple illustration of such relations and their association with the corresponding significance levels of tests. Note also that $b_{\cal T}(\Pi_{\theta_N}^N) \asymp \theta_N \sqrt N.$ 
\end{re}

\noindent
\begin{cor}\label{cor1} 
Let $m=\lfloor N\eta \rfloor,\; n=N-m,\; \Delta(N) = \lfloor N^p \rfloor,\; 0<p\leq 1/3,$ and $\theta_N=N^{-q}$ with $0<p<q<1/2$. Then the assumptions (i), (ii)' and (iii) of Theorem \ref{thm6} are satisfied for any $\nu < (q-p)/q$ while the corresponding significance levels, for a given $q$, satisfy $\alpha_N \asymp \exp\{-CN^{1-2q}\}$, where $C$ is a positive constant.
\end{cor}

\noindent
\begin{re}\label{re9}
The proofs of Theorems \ref{thm3} and \ref{thm5} indicate that $\Delta(N)$ only influences the results by defining the cut-off points from the ends of the interval (0,1). Theorem \ref{thm6} shows that we may consider an interval slightly narrower than $[1/\sqrt N,1-1/\sqrt N]$. The choice of the partition points, the $\pi_{jN}$'s, inside $[1/(\Delta(N) +1),
1-1/(\Delta(N)+1)]$ is practically immaterial to the final asymptotic results. In practice, it is natural to take a reasonably large number of such points to ensure accuracy, while ensuring that the calculations are not numerically complex. The line of our proofs also makes it clear that, instead of a discretized variant ${\cal W}_N$, one may simply consider
$$
\sup_{\epsilon(N) \leq t \leq 1-\epsilon(N)} \frac{(\hat G_n-\hat F_m)\circ \hat J_N^{-1}(t)}{\sqrt{t(1-t)}},
$$
where $\epsilon(N) \approx o(1/\sqrt N)$, and obtain similar asymptotic results.
\end{re}

\noindent
\section{ Simulation results}\label{S1}

\noindent
We compare the empirical powers of four tests which reject $\mathbb{H}_0$ for large values of the following statistics
$$
{\cal V}_N = \sqrt{\frac{mn}{N}} \sup_{z \in \mathbb{R}} \Bigl\{\hat G_n(z) - \hat F_m(z)\Bigr\},\;\;\;
$$
$$
{\cal T}_N^{\star} = \max_{1 \leq j \leq \Delta(N)} \bigl\{-{\cal L}_j\bigr\},\;\;\;
\mbox {with} \;\;\;\Delta(N)=2^{\lfloor \log_2 N \rfloor} -1,\;\;\pi_{jN}=\frac{j}{\Delta(N)+1},\;\;j=1,...,\Delta(N),
$$
$$
{\cal T}_N^o = \max_{1+\lfloor\sqrt N \rfloor \leq j \leq N-\lfloor\sqrt N \rfloor} \bigl\{-{\cal L}_j\bigr\},\;\;\;\mbox{where}\;\;\pi_{jN}=\frac{j}{N+1},\;\;j=1,...,N,
$$
and
$$
{\cal V}_N^e={\cal V}_{N^{(e)}},\;\;\mbox{where}\;N^{(e)}=m^{(e)}+n^{(e)}\;\,\mbox{and}\;m^{(e)}=\lfloor me_{{\cal T}{\cal V}}\rfloor,\;n^{(e)}=\lfloor ne_{{\cal T}{\cal V}}\rfloor.
$$
${\cal T}_N^{\star}$ is the variant of ${\cal T}_N$ which was studied in \citeauthor{r24} (\citeyear{r24,r26}), in Sections 3.2 and 2.2, respectively. ${\cal T}_N^o$ is a new test statistic which, in comparison to ${\cal T}_N^{\star}$, puts less weight on extreme observations. ${\cal V}_N^e$ is the two-sample Kolmogorov-Smirnov statistic based on samples of sizes which have been adapted to the intermediate efficiency value $e_{{\cal T}{\cal V}}=e_{{\cal T}{\cal V}}(\eta)$, as defined above. 

Taking into account the definition of the intermediate efficiency, it is expected that the empirical power of $ {\cal V}_N^e={\cal V}_{N^{(e)}}$ will be greater or equal to the corresponding powers of ${\cal T}_N^{\star}$ and 
${\cal T}_N^o $. We emphasize that in our simulation study, described below, the significance level and alternatives are independent of $N$. So, although the theory is local and assumes that significance levels tend to 0 as $N \to \infty$, the results demonstrate that this theory accurately describes the exact behavior of power under a standard simulation scheme. 

We consider the power and efficiency under some commonly used models of alternatives; cf. \cite{r4}, \cite{r11}, \cite{r23}, \cite{r36}.
 All of these models were considered in extensive simulation studies in \citeauthor{r24} (\citeyear{r24,r26}). This enables further comparisons. 

A description of the distributions used in our simulation study is given below. We list the models in the same order as they appear in Figures 1 - 4.\\
\begin{itemize}
\item $MU(a),\;\;a \in [-1,1],$ has density function $\;1 + a \sin(10 \pi x){\bf 1}(0.4 < x < 0.6),\; x \in [0,1]$; cf. Fan(1996).

\item $Pareto(a)$ coincides with $\;SM(a,1,1)$; see  below.

\item $SM(a,b,c),\;\;a \geq 1,\;b \geq 1,\;c \geq 1,\;$ denotes the Singh-Maddala model which obeys the distribution function $\;1-[1+(x/b)^a]^{-c},\;x > 0$.

\item $LN(a,b),\;a \in \mathbb{R},\; b>0,$ is the log-normal distribution with density function $\exp[-(\log x -a)^2/(2b^2)]/(\sqrt{2 \pi}bx),\;x>0,$ while $\beta LN(a_1,b_1)+(1-\beta)LN(a_2,b_2),\;\beta \in (0,1),$ denotes a mixture of such distributions.

\item $N(a,b)$ denotes the normal distribution with mean $a$ and dispersion $b$, $\chi_1^2$ denotes the central chi-square distribution with 1 degree of freedom, and $\beta N(a,b)+(1-\beta) \chi_1^2$ is a mixture of such distributions.

\item $Laplace(a,b),\;\;a \in \mathbb{R},\;b>0,\;$ has the density function $\;\exp(-|x-a|/b)/2b,\;x \in \mathbb{R}$.

\end{itemize}

The simulation results are presented graphically. The notation used in Figures 1 - 4 of this paper are identical to that previously used in the two papers by Ledwina and Wy{\l}upek mentioned above. In particular, given the alternative $(F_1,G_1)$, we indicate the underlying distributions in the two samples  by $F_1/G_1$. 
 
The number of Monte Carlo runs is 5\;000 throughout. The significance level was set to be $\alpha=0.01$ in all cases apart from the middle row in Figure 4, where it varies over the interval [0.001,0.010].

For each alternative $F_1/G_1$ under consideration, we plot in the first row of the corresponding figure 
the values of $e_{{\cal T}{\cal V}}=e_{{\cal T}{\cal V}}(\eta)$ (continuous brown line) and 
$\arg \max A^*(t;\eta)$ (blue dotted line) against $\eta \in (0,1)$. Theorem 6 shows that the location of $\arg \max A^*(t;\eta)$ on (0,1) is decisive to the magnitude of $e_{{\cal T}{\cal V}}(\eta)$.

Let us start with a description and some discussion of the results presented in Figures 1 - 3. For these figures, we took $N=800$ in the unbalanced case ($\eta_N \neq 0.5$) and various values of $N \in [200,500]$ for balanced  partitions ($\eta_N=0.5$). The above classification into balanced and unbalanced cases is not very precise, since among the unbalanced cases there is always a balanced one. However, such terminology allows a more succinct description of the results.

Figures 1 - 3 are ordered such that the maximal attainable value of $e_{{\cal T}{\cal V}}(\eta)$ over $\eta \in (0,1)$ is increasing. In Figures 1 - 3 the range of $\eta$ is  $[0.01 ,0.99]$. 
Observe that the maximal value of $e_{{\cal T}{\cal V}}(\eta)$ over the above mentioned interval ranges from 1.016, in the case of the alternative $MU(-1)/MU(1)$, to 22, for the  alternative $Laplace(0,1)/Laplace(1,1.25)$. Since the vertical scales in the first rows of Figures 1 - 3 are different in each case, in order to increase readability, we display the value of $e_{{\cal T}{\cal V}}(0.5)$ in the middle row of each figure.

In Figure 1, we illustrate the empirical powers of the four tests  for the selection of sample sizes described above, both for balanced and unbalanced partitions,  under two pairs of alternative distributions, corresponding to the Fan and Pareto models, respectively. The middle row contains the results for balanced partitions, while the bottom row describes the results for unbalanced ones. Figures 2 and 3 are constructed in an analogous manner. 

The behavior of the test based on ${\cal V}_N^e$ illustrates how the intuitive meaning of the efficiency measure manifests itself for finite sample sizes. 
It is expected that when $e_{{\cal T}{\cal V}} > 1$, the empirical powers of ${\cal V}_N^e$ should be slightly higher than the corresponding powers of ${\cal T}_N^{\star}$ and ${\cal T}_N^o$. We see that indeed this is the case regardless of the model, exact form of the efficiency function and whether the partition is balanced or not. 
We considered $\eta_N \in [0.1,0.9]$. For moderately large sample sizes, as used in the cases presented in Figures 1 - 3,  the accuracy of the prediction of the empirical power of ${\cal V}_N^{e}$ is very high for $\eta_N \in [0.2,0.8]$.

In the cases where $e_{{\cal T}{\cal V}}(\eta)$ is 1 or very close to 1 for all $\eta$ (cf. Figure 1), the empirical powers of ${\cal V}_N$ may be greater than the corresponding powers of ${\cal T}_N^{\star}$ and ${\cal T}_N^{o}$. The alternative that we consider in the first column of Figure 1, based on Fan (1996), may serve as an illustration of such a situation. In a sense, this is the least favorable situation for weighted statistics which are designed to be sensitive to differences between tails. Indeed, the two distributions $F_1$ and $G_1$ have the same tails and differ only in the central part. In such a case, the simpler  structure of ${\cal V}_N$ plays a role.
It can also be seen from Figure 1 that in a slightly less extreme case, represented by the pair of Pareto distributions, this deficiency disappears.

The evidence provided in Figures 1 - 3 clearly shows that, except in some very difficult circumstances (very small or very large $\eta_N$), applying the concept of the intermediate efficiency  to finite samples works very well. This situation should obviously be even better when $N$ is larger; cf. Figure  4 and the related comments given below.

We also studied larger sample sizes and smaller significance levels. In Figure 4 we display the results of such experiments for two models: one, corresponding to Singh-Maddala distributions,  where the efficiency is very close to 1 for all $\eta$, and the second one, corresponding to log-normal distributions, where the efficiency lies in an interval approximately 
equal to [4,16], depending on the value of $\eta$. In the middle row, we present empirical powers against $\alpha \in [0.001,0.010]$ for balanced partitions and $N=16\; 000$. We see that the efficiency gives very accurate and stable results. The same comment is valid for empirical powers under unbalanced partitions and $\alpha=0.01$. These results are presented in the bottom panels of Figure 4. It can be seen that for all $\eta \in [0.1,0.9]$ the simulation results are satisfactory.

To close, we comment on the differences between the empirical behavior of ${\cal T}_N^{\star}$ and ${\cal T}_N^o$.
Our simulations show that ${\cal T}_N^{\star}$ and ${\cal T}_N^o$ behave similarly for balanced partitions as long as the efficiency $e_{{\cal T}{\cal V}}$ is not very large. In the opposite case (cf. Figures 2 - 4), there is some gain from using ${\cal T}_N^{\star}$. 
 The explanation of this is simple. High efficiency occurs when arg max $A^{\star}(t;\eta)$ is located near 0 or 1. Obviously ${\cal T}_N^{\star}$ has a greater chance to detect such changes, as it allows closer inspection of weighted two-sample rank processes at arguments which are closer to 0 and 1 than ${\cal T}_N^{o}$ does; cf. ({3.9}) and related comments. In the cases of highly unbalanced partitions and small or moderate efficiencies, the new solution ${\cal T}_N^o$ provides some improvement under most of the alternatives.
It is also worth noticing that the value of $\Delta(N)$ corresponding to ${\cal T}_N^o$ approximately satisfies the requirements of our theoretical results. For ${\cal T}_N^{\star}$, the choice of $\Delta(N)$ is outside the allowable range. In any case,  applying the concept of intermediate efficiency to finite samples works well for ${\cal T}_N^{\star}$ in our simulation experiments.


\begin{figure}[H]
\centering
\includegraphics[scale=0.6]{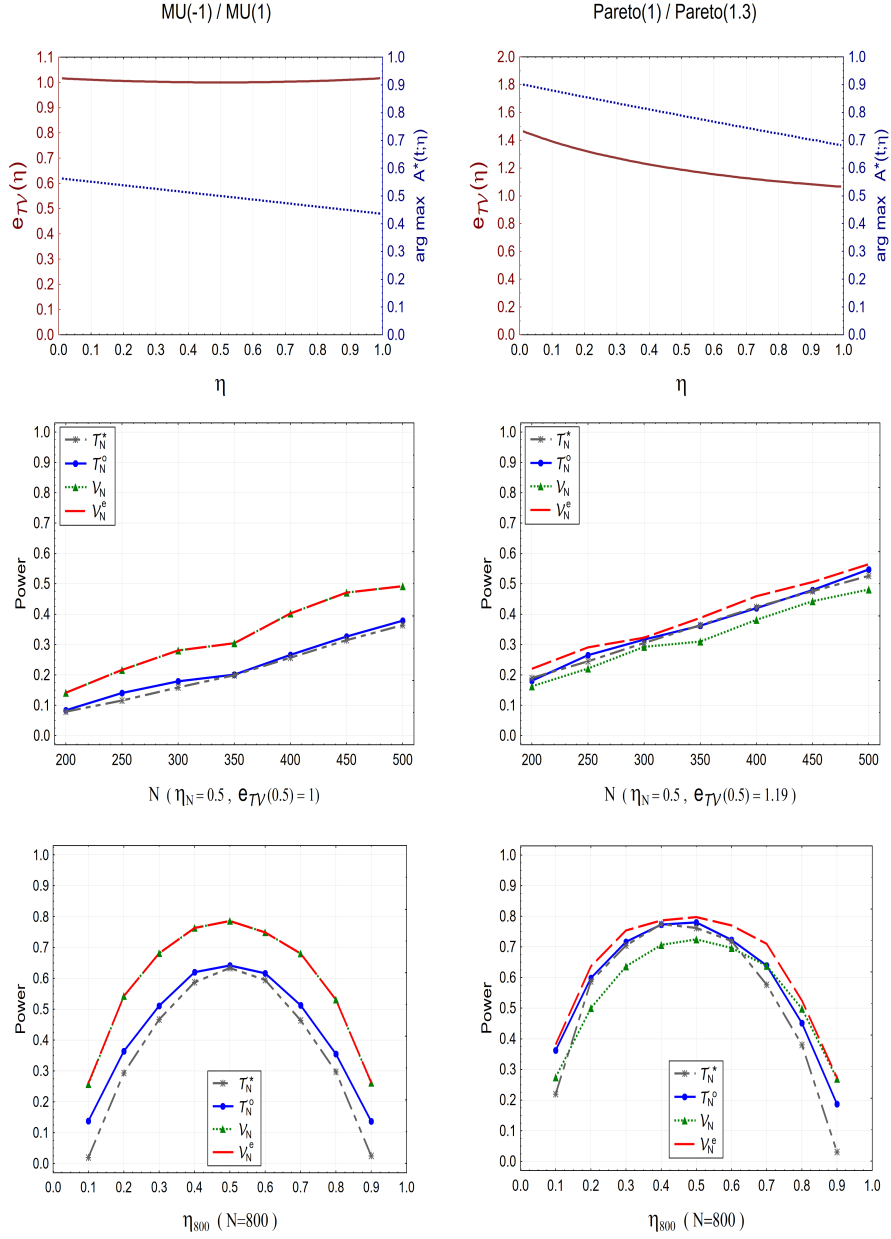} 
\vspace{-0.1cm}
\caption{{\bf  Fan and Pareto alternatives}. 
Moderately large sample sizes, $\alpha=0.01$.
{\it Upper panels}: efficiencies $e_{{\cal T}{\cal V}}(\eta)$ against $\eta \in (0,1)$ - brown continuous line; locations of the maximum of $A^*(t,\eta)$ against $\eta \in (0,1)$ - blue dotted line. {\it Middle panels}: empirical powers for balanced partitions ($\eta_N=0.5$) against $N \in [200,500]$. {\it Bottom panels}: empirical powers for unbalanced partitions against $\eta_N \in [0.1,0.9]$; $N=800$.}

\end{figure}

\begin{figure}[H]
\centering

\includegraphics[scale=0.6]{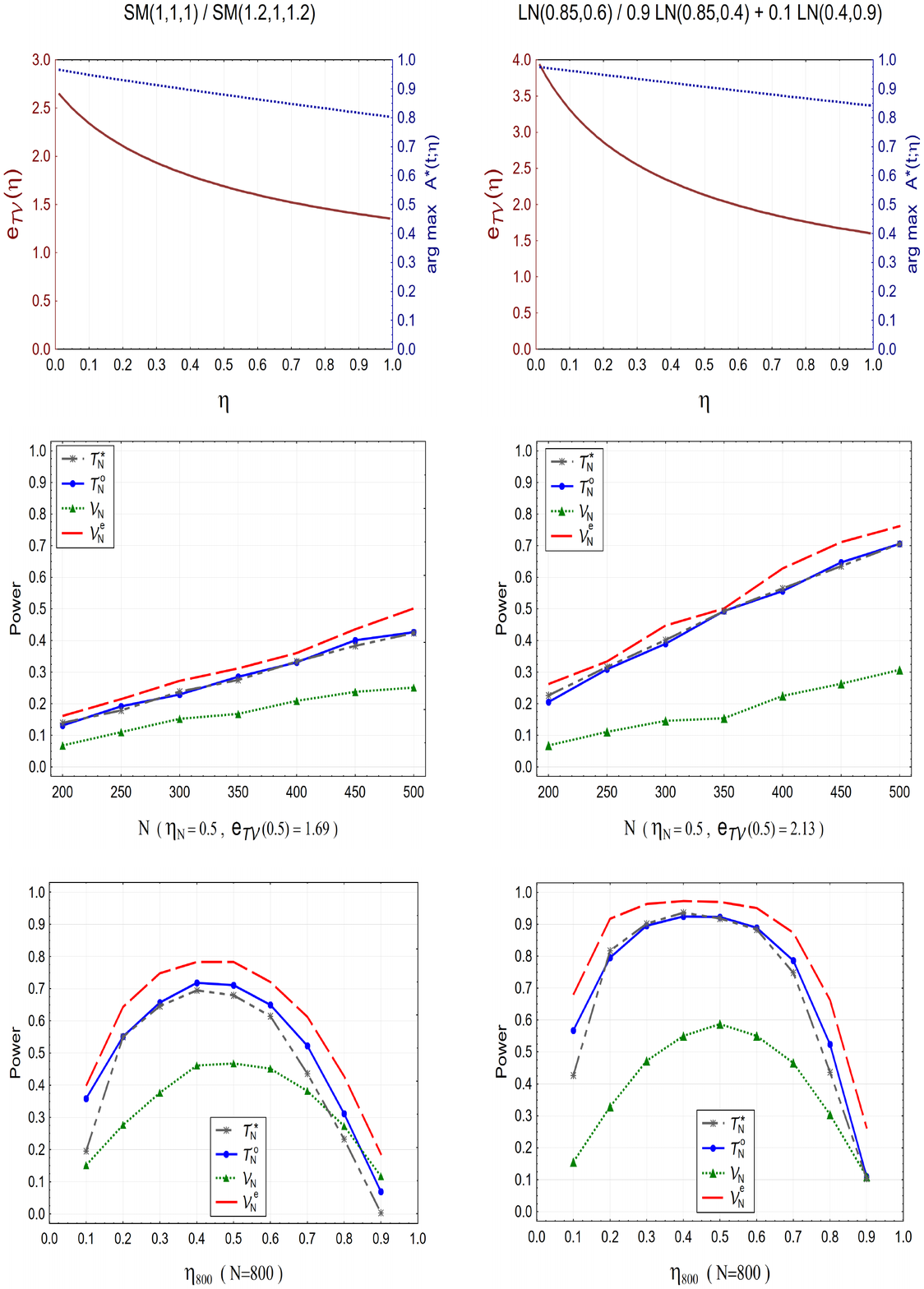}

\vspace{0.5cm}
\caption{{\bf  Singh-Maddala and log-normal alternatives}. Moderately large sample sizes, $\alpha=0.01$.
{\it Upper panels}: efficiencies $e_{{\cal T}{\cal V}}(\eta)$ against $\eta \in (0,1)$ - brown continuous line; locations of the maximum of $A^*(t,\eta)$ against $\eta \in (0,1)$ - blue dotted line. {\it Middle panels}: empirical powers for balanced partitions ($\eta_N=0.5$) against $N \in [200,500]$. {\it Bottom panels}: empirical powers for unbalanced partitions against $\eta_N \in [0.1,0.9]$; $N=800$.}

\end{figure}

\begin{figure}[H]
\centering

\includegraphics[scale=0.6]{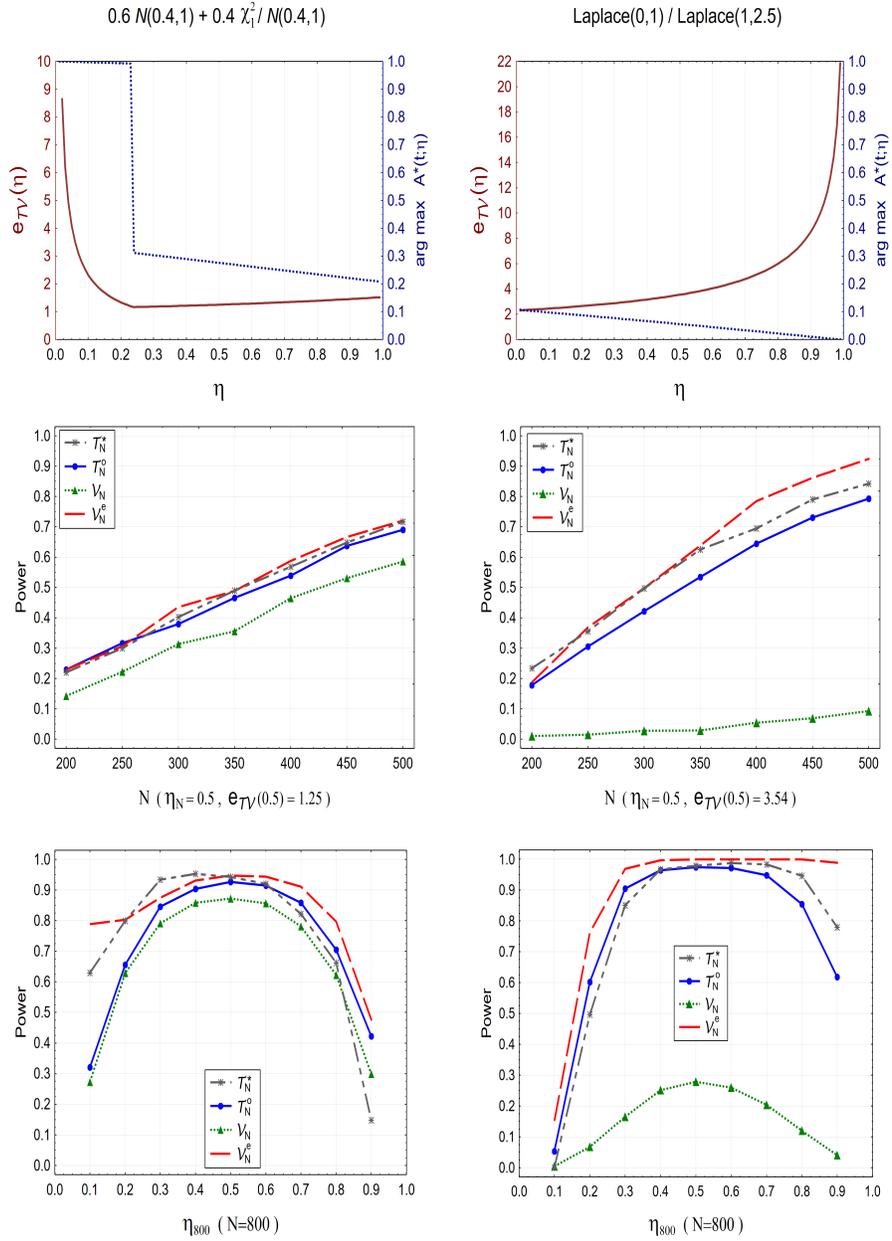}

\vspace{0.5cm}
\caption{{\bf Normal-chi-square mixture  and Laplace alternatives}.  Moderately large sample sizes, $\alpha=0.01$.
{\it Upper panels}: efficiencies $e_{{\cal T}{\cal V}}(\eta)$ against $\eta \in (0,1)$ - brown continuous line; locations of the maximum of $A^*(t,\eta)$ against $\eta \in (0,1)$ - blue dotted line. {\it Middle panels}: empirical powers for balanced partitions ($\eta_N=0.5$) against $N \in [200,500]$. {\it Bottom panels}: empirical powers for unbalanced partitions against $\eta_N \in [0.1,0.9]$; $N=800$.}

\end{figure}





\begin{figure}[H]
\centering
\includegraphics[scale=0.6]{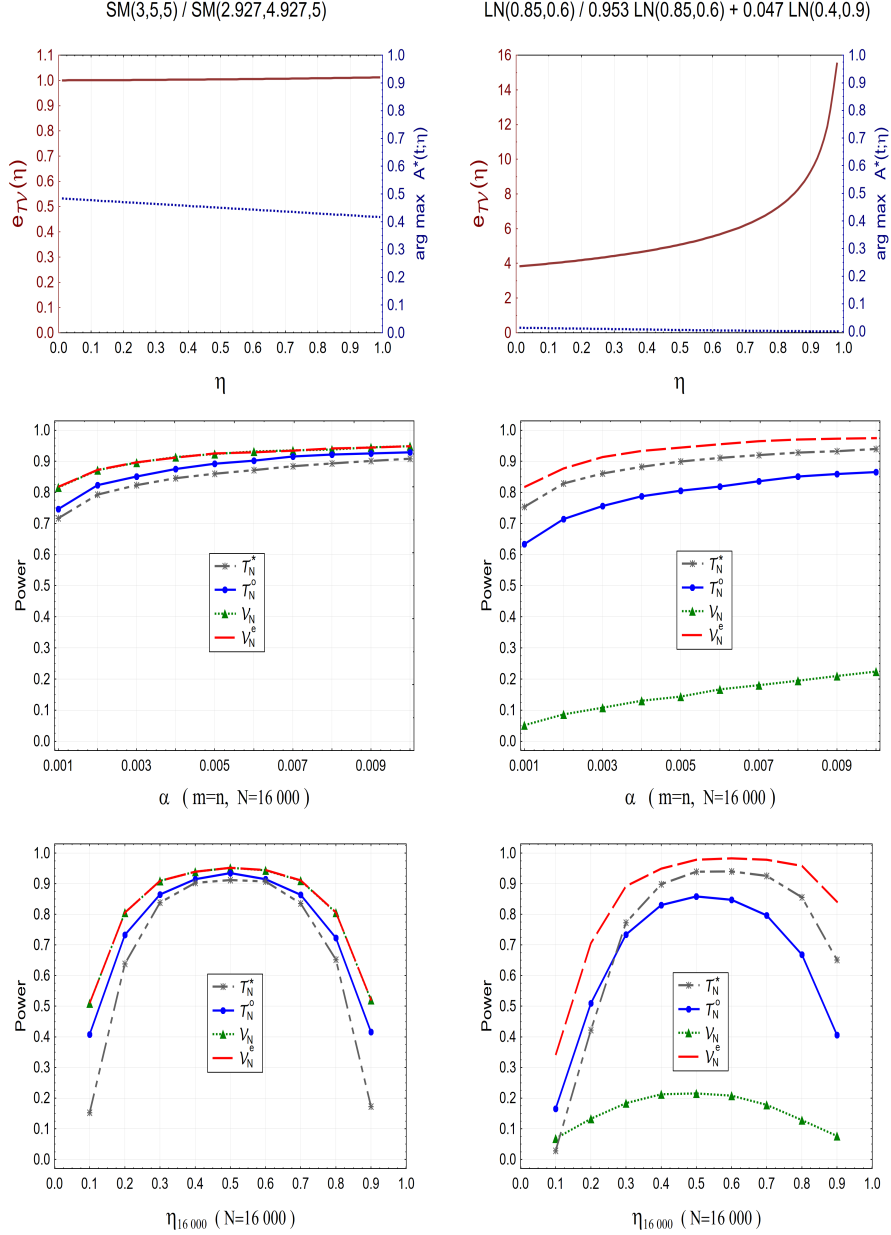} 
\vspace{-0.1cm}
\caption{{\bf Singh-Maddala and log-normal alternatives}. Large sample sizes. {\it Upper panels}: efficiencies $e_{{\cal T}{\cal V}}(\eta)$ against $\eta \in (0,1)$ - brown continuous line;  locations of the maximum of $A^*(t,\eta)$ against $\eta \in (0,1)$ - blue dotted line. {\it Middle panels}: empirical powers for balanced partitions ($\eta_N=0.5$) against $\alpha \in [0.001,0.010]$; $N=16\;000$. {\it Bottom panels}: empirical powers for unbalanced partitions against $\eta_N$; $N=16\;000$; $\alpha=0.01$.}

\end{figure}

\begin{center}
APPENDICES
\end{center}

\noindent{\bf Appendix A: Proof of Theorem 1}\\

\noindent
Let $\{\alpha_N\}$ be any sequence from $\mathbb{L}^*$ satisfying (2.11).\\

\noindent
{\bf Step 1. Basic relation between $\alpha_N$, $c_{{\cal U}^{(II)}}$ and $b
_{{\cal U}^{(II)}}(\cdot)$.} 
For any $h>0$ define $w_N=\sqrt{(-h\log \alpha_N)/Nc_{{\cal U}^{(II)}}}$. Since (2.11) holds therefore for this $w_N$ (2.9) applies and yields 
$$
\frac{c_{{\cal U}^{(II)}}}{h\log\alpha_N}\log\sup_{P\times Q\in\mathbb{P}_0} P^{m(N)}\times Q^{n(N)}\bigl({\cal U}_N^{(II)}\geq w_N\sqrt{N}\bigr)\to c_{{\cal U}^{(II)}}. 
\eqno (A.1)
$$
For arbitrary $\epsilon>0$ take $h=1+\epsilon$ in $w_N$ appearing in (A.1). Then, for sufficiently large $N$,
$$ 
\sup_{P\times Q\in\mathbb{P}_0} P^{m(N)}\times Q^{n(N)}\bigl({\cal U}_N^{(II)}\geq w_N\sqrt{N}\bigr)\leq \alpha_N
$$ 
which, by the definition of $u_{\alpha N}^{(II)}$, means that $u_{\alpha_N N}^{(II)}\leq \sqrt{(-(1+\epsilon)\log \alpha_N)/c_{{\cal U}^{(II)}}}$ or equivalently $c_{{\cal U}^{(II)}}[u_{\alpha_NN}^{(II)}]^2\leq -(1+\epsilon)\log\alpha_N$.
Similarly taking in (A.1) $h=1-\epsilon$ we get $c_{{\cal U}^{(II)}}[u_{\alpha_NN}^{(II)}]^2\geq -(1-\epsilon)\log\alpha_N$ for sufficiently large $N$. Since $\epsilon$ was taken arbitrarily we obtain
$$ 
-\log\alpha_N=c_{{\cal U}^{(II)}} [u_{\alpha_NN}^{(II)}]^2[1+o(1)]. 
\eqno (A.2)
$$

On the other hand, by (2.10), we have for arbitrary $\epsilon>0$
$$ 
\Pi_{\theta_N}^N\bigl({\cal U}_N^{(II)}\geq (1+\epsilon)b_{{\cal U}^{(II)}}(\Pi_{\theta_N}^N)\bigr)\to 0, 
\;\;\Pi_{\theta_N}^N\bigl({\cal U}_N^{(II)}\geq (1-\epsilon)b_{{\cal U}^{(II)}}(\Pi_{\theta_N}^N)\bigr)\to 1.
$$
Since $\{\alpha_N\}\in \mathbb{L}^*$, the condition (2.4) implies $(1-\epsilon)b_{{\cal U}^{(II)}}(\Pi_{\theta_N}^N)\leq u_{\alpha_NN}^{(II)}\leq (1+\epsilon)b_{{\cal U}^{(II)}}(\Pi_{\theta_N}^N)$ for sufficiently large $N$. As $\epsilon$ is arbitrary this means that $u_{\alpha_NN}^{(II)}=b_{{\cal U}^{(II)}}(\Pi_{\theta_N}^N)[1+o(1)]$ which together with (A.2) gives
$$
-\log\alpha_N=c_{{\cal U}^{(II)}}\bigl[b_{{\cal U}^{(II)}}(\Pi_{\theta_N}^N)\bigr]^2[1+o(1)]. 
\eqno (A.3)
$$
From (2.12) and (A.3) we have
$$
-\log\alpha_N={\bf e} c_{{\cal U}^{(I)}}\frac{m(N)n(N)}{N}\bigl[b_{{\cal U}^{(I)}}(P_{\theta_N}\times Q_{\theta_N})\bigr]^2[1+o(1)]
\;\;\;\mbox{if}\;\;{\bf e}\in (0,\infty),
\eqno (A.4)
$$
$$-\log\alpha_N=o\left(\frac{mn}{N}\bigl[b_{{\cal U}^{(I)}}(P_{\theta_N}\times Q_{\theta_N})\bigr]^2\right) \;\;\;\mbox{if}\;\;{\bf e}=0,
\eqno (A.5)
$$
and
$$\frac{mn}{N}\bigl[b_{{\cal U}^{(I)}}(P_{\theta_N}\times Q_{\theta_N})\bigr]^2=o(-\log\alpha_N) \;\;\;\mbox{if}\;\;{\bf e}=\infty.
\eqno (A.6)
$$
\\

\noindent
{\bf Step 2. Lower bound for the fraction of sample sizes.}
For ${\bf e} \in (0,\infty]$ we shall show that (cf. (2.7))
$$ 
\liminf_{N\to\infty}\frac{M_{{\cal U}^{(II)}{\cal U}^{(I)}}\bigl(N,\Scale[1.6]{\pi}\bigr)}{N}\geq {\bf e}. 
\eqno (A.7)
$$
Suppose, contrary, that there exists an increasing sequence $\{k_j\}$ of natural numbers such that $k_j \to \infty$ as $j \to \infty$ and 
$$
\frac{M_{{\cal U}^{(II)}{\cal U}^{(I)}}\bigl(k_j,\Scale[1.6]{\pi}\bigr)}{k_j}\to 
\gamma < {\bf e}.
$$
For  $\tau>0$ such that $\gamma+\tau<{\bf e}$,  define
$$
\gamma_j=\frac{M_{{\cal U}^{(II)}{\cal U}^{(I)}}\bigl(k_j,\Scale[1.6]{\pi}\bigr)}{k_j}+\frac{\lfloor\tau k_j\rfloor}{k_j} .
$$
 Then $\lim_{j \to \infty} \gamma_j = \gamma +\tau \in (0, {\bf e})$.
Moreover, $\{\gamma_jk_j\}$ is a sequence of integers and $\gamma_jk_j\geq M_{{\cal U}^{(II)}{\cal U}^{(I)}}(k_j,\Scale[1.6]{\pi})$ for sufficiently large $j$. Hence by (2.6)
$$ 
\Pi_{\theta_{k_j}}^{\gamma_j k_j}\bigl({\cal U}_{\gamma_jk_j}^{(I)}>u_{\alpha_{k_j}\gamma_jk_j}^{(I)}\bigr)\geq 
\Pi_{\theta_{k_j}}^{k_j}\bigl({\cal U}_{k_j}^{(II)}>u_{\alpha_{k_j}k_j}^{(II)}\bigr). 
\eqno (A.8)
$$
Since $\{\gamma_j\}$ has positive and finite limit and $\{\alpha_N\}\in\mathbb{L}$ the assumption (I.1) can be applied to the sequence $\{w_N\}$ defined as follows: $w^2_{\gamma_jk_j}=-(1-\delta)(\log \alpha_{k_j})/\gamma_jk_jc_{{\cal U}^{(I)}}$ for $j=1,2,...$ and $w_N^2=-(1-\delta)(\log\alpha_N)/Nc_{{\cal U}^{(I)}}$ for $N\neq \gamma_jk_j$, where $\delta \in (0,1)$ is arbitrary. 
The assumption (I.1) applied for the subsequence $\{\gamma_jk_j\}$ yields
$$
\frac{c_{{\cal U}^{(I)}}}{(1-\delta)\log \alpha_{k_j}}\log\sup_{P\times Q\in\mathbb{P}_0} P^{m(\gamma_jk_j)}\times Q^{n(\gamma_jk_j)}\left({\cal U}_{\gamma_jk_j}^{(I)}\geq \sqrt{\frac{-(1-\delta)\log\alpha_{k_j}}{c_{{\cal U}^{(I)}}}}\right)\to c_{{\cal U}^{(I)}},
$$ 
which means that for sufficiently large $j$
$$
\sup_{P\times Q\in\mathbb{P}_0} P^{m(\gamma_jk_j)}\times Q^{n(\gamma_jk_j)}\left({\cal U}_{\gamma_jk_j}^{(I)}\geq \sqrt{\frac{-(1-\delta)\log\alpha_{k_j}}{c_{{\cal U}^{(I)}}}}\right)\geq \alpha_{k_j}.
$$
This and the definition of $u_{\alpha N}^{(I)}$ imply
$$u_{\alpha_{k_j}\gamma_jk_j}^{(I)}\geq \sqrt{\frac{-(1-\delta)\log\alpha_{k_j}}{c_{{\cal U}^{(I)}}}}.$$
Hence, from (A.8) and the fact that $\{\alpha_N\}\in \mathbb{L}^*$ we obtain
$$ 
\liminf_{j\to\infty}\Pi_{\theta_{k_j}}^{\gamma_jk_j}\left({\cal U}_{\gamma_jk_j}^{(I)}\geq \sqrt{\frac{-(1-\delta)\log\alpha_{k_j}}{c_{{\cal U}^{(I)}}}}\right) >0. 
\eqno (A.9)
$$

Now, consider a sequence $\{\vartheta_N\}$, being a modification of $\{\theta_N\}$, and defined as follows: $\vartheta_{\gamma_jk_j}=\theta_{k_j}$ for $j=1,2,...$ and $\vartheta_N=\theta_N$ for $N\neq \gamma_jk_j$. 

Since $\gamma_j\to\gamma+\tau\in(0,{\bf e})$ we have $\vartheta_N\to 0$ and $N\vartheta_N^{\rho}\to\infty$ and (I.2) can be used for this sequence. 
Hence (2.8) applied  to the subsequence $\{\gamma_jk_j\}$ and arbitrary $\epsilon>0$ implies
$$ 
\Pi_{\theta_{k_j}}^{\gamma_jk_j}\left({\cal U}_{\gamma_jk_j}^{(I)}\geq 
(1+\epsilon)\sqrt{\frac{m(\gamma_jk_j)n(\gamma_jk_j)}{\gamma_jk_j}}b_{{\cal U}^{(I)}}(P_{\theta_{k_j}}\times Q_{\theta_{k_j}})\right)\to 0.
\eqno (A.10)
$$

We shall show that the relations (A.9) and (A.10) give a contradiction.

Assume first that ${\bf e} \in (0,\infty)$. 
We have $\gamma+\tau<{\bf e}$. For fixed $\delta<1-(\gamma+\tau)/{\bf e}$ choose $\epsilon>0$ so small that $(1+\epsilon)^2(\gamma+\tau)<(1-\delta){\bf e}$. 
Define $\kappa_N = \sqrt{\eta_N(1-\eta_N)}=\sqrt{m(N)n(N)}/N$.
By (A.4), the convergence $\kappa_N\to\kappa$ and $\gamma_j\to\gamma+\tau$, and the choice of $\delta$ and $\epsilon$ we have, for $j$ sufficiently large,
$$ 
\frac{-(1-\delta)\log \alpha_{k_j}}{c_{{\cal U}^{(I)}}}=(1-\delta){\bf e}k_j\kappa_{k_j}^2\bigl[b_{{\cal U}^{(I)}}(P_{\theta_{k_j}}\times Q_{\theta_{k_j}})\bigr]^2[1+o(1)]\hspace{3cm}$$
$$\hspace{3.5cm}=(1-\delta)\frac{{\bf e}}{\gamma+\tau} \frac{\kappa_{k_j}^2}{\kappa_{\gamma_jk_j}^2}\frac{m(\gamma_jk_j)n(\gamma_jk_j)}{\gamma_jk_j}\bigl[b_{{\cal U}^{(I)}}(P_{\theta_{k_j}}\times Q_{\theta_{k_j}})\bigr]^2[1+o(1)]$$
$$>(1+\epsilon)^2\frac{m(\gamma_jk_j)n(\gamma_jk_j)}{\gamma_jk_j}\bigl[b_{{\cal U}^{(I)}}(P_{\theta_{k_j}}\times Q_{\theta_{k_j}})\bigr]^2\hspace{0.2cm}$$
which contradicts (A.9) and (A.10).

If ${\bf e}=\infty$ we have from (A.6) and the convergence $\gamma_j\to\gamma+\tau$ and $\kappa_N\to\kappa=\sqrt{\eta(1-\eta)}$ 
$$
\frac{m(\gamma_jk_j)n(\gamma_jk_j)}{\gamma_jk_j}\bigl[b_{{\cal U}^{(I)}}(P_{\theta_{k_j}}\times Q_{\theta_{k_j}})\bigr]^2=\frac{\kappa^2_{\gamma_jk_j}}{\kappa^2_{k_j}}\gamma_j\frac{m(k_j)n(k_j)}{k_j}\bigl[b_{{\cal U}^{(I)}}(P_{\theta_{k_j}}\times Q_{\theta_{k_j}})\bigr]^2
$$
$$
=o\Bigl(\frac{-(1-\delta)\log \alpha_{k_j}}{c_{{\cal U}^{(I)}}}\Bigr)$$
which contradicts (A.9) and (A.10), as well.\\

\noindent
{\bf Step 3. Upper bound for the fraction of sample sizes.}
For ${\bf e}\in [0,\infty)$ we shall show that
$$ 
\limsup_{N\to\infty}\frac{M_{{\cal U}^{(II)}{\cal U}^{(I)}}\bigl(N,\Scale[1.6]{\pi}\bigr)}{N}\leq {\bf e}. 
\eqno (A.11)
$$
The argument is very similar to that of Step 2. 
Suppose,  that there exists an increasing sequence $\{k_j\}$ of natural numbers such that
$$
\gamma_j=\frac{M_{{\cal U}^{(II)}{\cal U}^{(I)}}\bigl(k_j,\Scale[1.6]{\pi}\bigr)-1}{k_j}\to \gamma > {\bf e}.
$$
Note that $\gamma$ may be equal to $\infty$. Since $\gamma_jk_j=M_{{\cal U}^{(II)}{\cal U}^{(I)}}(k_j,\Scale[1.6]{\pi})-1$ then by (2.6)
$$ 
\Pi_{\theta_{k_j}}^{\gamma_jk_j}\bigl({\cal U}_{\gamma_jk_j}^{(I)}>u_{\alpha_{k_j}\gamma_jk_j}^{(I)}\bigr)< \Pi_{\theta_{k_j}}^{k_j}\bigl({\cal U}_{\gamma_jk_j}^{(II)}>u_{\alpha_{k_j}k_j}^{(II)}\bigr). 
\eqno (A.12)
$$
Since $\{\gamma_j\}$ has positive limit or tends to $\infty$ and $\{\alpha_N\}\in\mathbb{L}$, then the condition (I.1) can be applied to the sequence $\{w_N\}$ defined as follows: $w^2_{\gamma_jk_j}=-(1+\delta)(\log \alpha_{k_j})/\gamma_jk_jc_{{\cal U}^{(I)}}$ for $j=1,2,...$ and $w_N^2=-(1+\delta)(\log\alpha_N)/Nc_{{\cal U}^{(I)}}$ for $N\neq \gamma_jk_j$, where $\delta>0$ is arbitrary. By (I.1) applied to the subsequence $\{\gamma_jk_j\}$ we get
$$
\frac{c_{{\cal U}^{(I)}}}{(1+\delta)\log \alpha_{k_j}}\log\sup_{P\times Q\in\mathbb{P}_0} P^{m(\gamma_jk_j)}\times Q^{n(\gamma_jk_j)}\left({\cal U}_{\gamma_jk_j}^{(I)}\geq \sqrt{\frac{-(1+\delta)\log\alpha_{k_j}}{c_{{\cal U}^{(I)}}}}\right)\to c_{{\cal U}^{(I)}}
$$ 
which means that for sufficiently large $j$
$$
\sup_{P\times Q\in\mathbb{P}_0} P^{m(\gamma_jk_j)}\times Q^{n(\gamma_jk_j)}\left({\cal U}_{\gamma_jk_j}^{(I)}\geq \sqrt{\frac{-(1+\delta)\log\alpha_{k_j}}{c_{{\cal U}^{(I)}}}}\right)\leq \alpha_{k_j}.
$$
This and the definition of $u_{\alpha N}^{(I)}$ implies
$$
u_{\alpha_{k_j}\gamma_jk_j}^{(I)}\leq \sqrt{\frac{-(1+\delta)\log\alpha_{k_j}}{c_{{\cal U}^{(I)}}}}.
$$
Hence, from (A.12) and the fact that $\{\alpha_N\}\in \mathbb{L}^*$ we obtain
$$ 
\limsup_{j\to\infty}\Pi_{\theta_{k_j}}^{\gamma_jk_j}\left({\cal U}_{\gamma_jk_j}^{(I)}\geq \sqrt{\frac{-(1+\delta)\log\alpha_{k_j}}{c_{{\cal U}^{(I)}}}}\right)<1. 
\eqno (A.13)
$$

Now, consider a sequence $\{\vartheta_N\}$, being a modification of $\{\theta_N\}$, which is defined as follows: $\vartheta_{\gamma_jk_j}=\theta_{k_j}$ for $j=1,2,...$ and $\vartheta_N=\theta_N$ for $N\neq \gamma_jk_j$, where $\{k_j\}$ is the sequence selected at the beginning of this step. 

Since $\gamma\in({\bf e},\infty]$ we have $\vartheta_N\to 0$  and $N\vartheta_N^{\rho}\to\infty$ and (I.2) can be applied to this sequence. Hence (2.8) applied for the subsequence $\{\gamma_jk_j\}$ and arbitrary $\epsilon>0$ yields
$$ 
\Pi_{\theta_{k_j}}^{\gamma_jk_j}\left({\cal U}_{\gamma_jk_j}^{(I)}\geq (1-\epsilon)\sqrt{\frac{m(\gamma_jk_j)n(\gamma_jk_j)}{\gamma_jk_j}}b_{{\cal U}^{(I)}}(P_{\theta_{k_j}}\times Q_{\theta_{k_j}})\right)\to 1. 
\eqno (A.14)
$$

We shall argue that the relations (A.13) and (A.14) give a contradiction.

Indeed, if ${\bf e} \in (0,\infty)$ we have $\gamma>{\bf e}$. For fixed $\delta<\gamma/{\bf e}-1$ choose $\epsilon>0$ so small that $(1-\epsilon)^2\gamma>(1+\delta){\bf e}$. This, (A.4), $\kappa_N\to\kappa$ and $\gamma_j\to\gamma$ imply for $j$ sufficiently large 
$$ 
\frac{-(1+\delta)\log \alpha_{k_j}}{c_{{\cal U}^{(I)}}}=(1+\delta){\bf e}k_j\kappa_{k_j}^2[b_{{\cal U}^{(I)}}(P_{\theta_{k_j}}\times Q_{\theta_{k_j}})]^2[1+o(1)]\hspace{3cm}$$
$$\hspace{3cm}=(1+\delta)\frac{{\bf e}}{\gamma} \frac{\kappa_{k_j}^2}{\kappa_{\gamma_jk_j}^2}\frac{m(\gamma_jk_j)n(\gamma_jk_j)}{\gamma_jk_j}\bigl[b_{{\cal U}^{(I)}}(P_{\theta_{k_j}}\times Q_{\theta_{k_j}})\bigr]^2[1+o(1)]$$
$$<(1-\epsilon)^2\frac{m(\gamma_jk_j)n(\gamma_jk_j)}{\gamma_jk_j}\bigl[b_{{\cal U}^{(I)}}(P_{\theta_{k_j}}\times Q_{\theta_{k_j}})\bigr]^2$$
which contradicts (A.13) and (A.14). 

If ${\bf e}=0$ we have from (A.5) and the convergence $\gamma_j\to\gamma>0$ and $\kappa_N\to\kappa$ 
$$ 
\frac{-(1+\delta)\log \alpha_{k_j}}{c_{{\cal U}^{(I)}}}= o\bigl(k_j\bigl[b_{{\cal U}^{(I)}}(P_{\theta_{k_j}}\times Q_{\theta_{k_j}})\bigr]^2\bigr)
$$
$$=o\left(\frac{1}{\kappa^2_{\gamma_jk_j}}\frac{1}{\gamma_j}\frac{m(\gamma_jk_j)n(\gamma_jk_j)}{\gamma_jk_j}\bigl[b_{{\cal U}^{(I)}}(P_{\theta_{k_j}}\times Q_{\theta_{k_j}})\bigr]^2\right)$$
which contradicts (A.13) and (A.14), as well. The proof is complete. \hfill $\Box$\\ 


\noindent
{\bf Appendix B: Lemma 1, Lemma 2, and proof of Lemma 1} \\

\noindent
{\bf Lemma 1.}
{\it Let $\{P_{\theta_N} \times Q_{\theta_N}\}$ be the particular sequence of alternatives under consideration. Suppose that there exist cumulative distribution functions $U_1^{(II)}$ and $U_2^{(II)}$ and positive sequences $\{a_N^{(II)}\}$ and $\{b_N^{(II)}\}$ such that $b_N^{(II)} \to \infty,
\;b_N^{(II)}/a_N^{(II)} \to \infty$ and
for some $w_0^{(II)} \in \mathbb{R}$  we have

\noindent
(i) $ \limsup_{N \to \infty} \Pi_{\theta_N}^N \Bigl(\frac{{\cal U}_N^{(II)}-b_N^{(II)} }{a_N^{(II)}} \leq w\Bigr) \leq U_2^{(II)}(w)\;\;$ for all $\;w \in \mathbb{R}$,

\noindent
(ii) $ \liminf_{N \to \infty} \Pi_{\theta_N}^N \Bigl(\frac{{\cal U}_N^{(II)}-b_N^{(II)} }{a_N^{(II)}} \leq w\Bigr) \geq U_1^{(II)}(w)\;\;$ for all $\;w \in [w_0^{(II)},\infty)$.
Then (II.2)  holds true with $b_{{\cal U}^{(II)} }(\Pi_{\theta_N}^N)=b_N^{(II)}$.\\
\noindent
Further suppose that $w_0^{(II)}$ is such that for some $w_1^{(II)} > w_0^{(II)}$ satisfying $0 < U_1^{(II)}(w_0^{(II)}) < U_2^{(II)}(w_1^{(II)}) <1$  it holds for $N$ sufficiently large \\
\noindent
(iii) $\displaystyle
\sup_{P \times Q \in \mathbb{P}_0} P^{m(N)} \times Q^{n(N)} \Bigl(\frac{{\cal U}_N^{(II)}-b_N^{(II)} }{a_N^{(II)}} > w_1^{(II)}\Bigr) $
$$< \sup_{P \times Q \in \mathbb{P}_0} P^{m(N)} \times Q^{n(N)} \Bigl(\frac{{\cal U}_N^{(II)}-b_N^{(II)} }{a_N^{(II)}} > w_0^{(II)}\Bigr).
$$

Finally assume that (II.1) holds for ${\cal U}_N^{(II)}$ and  for the above $\{b_N^{(II)}\}$ and $\{\gamma_N\}$, $\{\lambda_N\}$ appearing in (II.1) it holds $[b_N^{(II)}]^2/\lambda_N \to 0$ and $[b_N^{(II)}]^2/\gamma_N \to \infty$. Then, (2.11) is satisfied  with 
$$
\alpha_N = \sup_{P \times Q \in \mathbb{P}_0} P^{m(N)} \times Q^{n(N)} \Bigl(\frac{{\cal U}_N^{(II)}-b_N^{(II)} }{a_N^{(II)}} > w_1^{(II)}\Bigr)
\eqno (B.1) $$
while the asymptotic power of the pertaining test based on ${\cal U}^{(II)}$ lies in the interval $[1-U_2^{(II)}(w_1^{(II)}),1-U_1^{(II)}(w_0^{(II)})] $.}\\

For completeness we state also a simple analogue of Lemma 1 which may be useful for checking (I.2) for ${\cal U}_N^{(I)}$. Its proof is quite similar to that of Lemma 1, so we omit it.\\

\noindent
{\bf Lemma 2.}
{\it Let $\{P_{\vartheta_N} \times Q_{\vartheta_N}\}$ be arbitrary sequence of alternatives for which $\vartheta_N \to 0$ and $N\vartheta_N^{\rho} \to \infty,\;\rho \in [1,2]$. 
Set $b_N^{(I)}=\sqrt{mn/N}b_{{\cal U}^{(I)}}(P_{\vartheta_N}\times Q_{\vartheta_N})$, where the positive function $b_{{\cal U}^{(I)}}(P_{\vartheta}\times Q_{\vartheta})$ is defined for all $\vartheta \in (0,1)$.
Suppose that, for each $\{\vartheta_N\}$ as above, there exist  cumulative distribution functions $U_1^{(I)}$ and $U_2^{(I)}$ and a positive sequence $\{a_N^{(I)}\}$ such that  $\;b_N^{(I)}/a_N^{(I)} \to \infty$ and for some $w_0^{(I)} \in \mathbb{R}$  we have

\noindent
(i) $ \limsup_{N \to \infty} \Pi_{\vartheta_N}^N \Bigl(\frac{{\cal U}_N^{(I)}-b_N^{(I)} }{a_N^{(I)}} \leq w\Bigr) \leq U_2^{(I)}(w)\;\;$ for all $\;w \in \mathbb{R}$,

\noindent
(ii) $ \liminf_{N \to \infty} \Pi_{\vartheta_N}^N \Bigl(\frac{{\cal U}_N^{(I)}-b_N^{(I)} }{a_N^{(I)}} \leq w\Bigr) \geq U_1^{(I)}(w)\;\;$ for all $\;w \in [w_0^{(I)},\infty)$.
\\
Then (I.2) is satisfied with the above $b_{{\cal U}^{(I)}}(\cdot)$. Distribution functions $U_1^{(I)}$ and $U_2^{(I)}$, the sequence $\{a_N^{(I)}\}$  as well as $w_0^{(I)}$ may be different for each sequence $\{\vartheta_N\}$.}\\

Note that if for all $N$ sufficiently large $\;\sup_{P \times Q \in \mathbb{P}_0} P^{m(N)} \times Q^{n(N)} \bigl({\cal U}_N^{(II)}-b_N^{(II)}  < w a_N^{(II)} \bigr)$ is strictly increasing in $w$  then (iii) holds true for every $w_1^{(II)} > w_0^{(II)}$.\\

\noindent
{\bf Proof of Lemma 1.}
First we shall check that indeed the conditions (i) and (ii) yield (2.10) with $b_{{\cal U}^{(II)}}(\Pi_{\theta_N}^N)=b_N^{(II)}$. We have
$$
\Pi_{\theta_N}^N\Bigl(\Bigl|\frac{{\cal U}_N^{(II)}}{b_N^{(II)}}-1\Bigr| \geq \epsilon\Bigr)=
\Pi_{\theta_N}^N\Bigl(\Bigl|\frac{{\cal U}_N^{(II)}-{b_N^{(II)}}}{a_N^{(II)}}\Bigr| \geq \frac{b_N^{(II)}}{a_N^{(II)}}\epsilon\Bigr).
\eqno(B.2)
$$
Since $b_N^{(II)}/a_N^{(II)} \to \infty$ we can take $w^* > w_0^{(II)}$, and $N$ enough large,  to majorize limsup of (B.2) by $1-U_1^{(II)}(w^*)+U_2^{(II)}(-w^*)$. Since $w^*$ can be arbitrary large the bound is arbitrary small. 

Now we shall check that $\{\alpha_N\}$ given in (B.1) fulfills the requirements needed to calculate the intermediate efficiency via Theorem 1.

By the definition of $u_{\alpha_{N}N}^{(II)}$ and (iii) it follows 
$$
b_N^{(II)} + a_N^{(II)} w_0^{(II)} \leq u_{\alpha_{N}N}^{(II)} \leq b_N^{(II)} + a_N^{(II)} w_1^{(II)}.
\eqno(B.3)
$$
Due to the assumptions on the sequences $\{a_N^{(II)}\}$, $\{b_N^{(II)}\}$, $\{\gamma_N\}$ and $\{\lambda_N\}$ we have
$$
(b_N^{(II)}+a_N^{(II)} w_1^{(II)})^2=[b_N^{(II)}]^2\Bigl(1+\frac{a_N^{(II)}}{b_N^{(II)}}w_1^{(II)}\Bigr)^2 \to \infty
$$
and
$$
\frac{(b_N^{(II)}+a_N^{(II)} w_1^{(II)})^2}{\lambda_N}=\frac{[b_N^{(II)}]^2}{\lambda_N}\Bigl(1+\frac{a_N^{(II)}}{b_N^{(II)}}w_1^{(II)}\Bigr)^2 \to 0,
$$
$$
\frac{(b_N^{(II)}+a_N^{(II)} w_1^{(II)})^2}{\gamma_N}=\frac{[b_N^{(II)}]^2}{\gamma_N}\Bigl(1+\frac{a_N^{(II)}}{b_N^{(II)}}w_1^{(II)}\Bigr)^2 \to \infty.
$$
Hence, for $w_N^2=(b_N^{(II)}+a_N^{(II)} w_1^{(II)})^2/N$ the condition (II.1) can be applied and yields
$$
-\frac{1}{(b_N^{(II)}+a_N^{(II)} w_1^{(II)})^2} \log \sup_{{P \times Q} \in \mathbb{P}_0} P^{m(N)} \times Q^{n(N)} \bigl({\cal U}_N^{(II)} \geq b_N^{(II)} + a_N^{(II)} w_1^{(II)}\bigr) \to c_{{\cal U}^{(II)}}.
$$
By (B.3) and the definition of $\alpha_N$, this implies
$$
-\frac{\log \alpha_N}{(b_N^{(II)}+a_N^{(II)} w_1^{(II)})^2} \to c_{{\cal U}^{(II)}}. 
$$
Similar argument works for $w_1^{(II)}$ replaced by $w_0^{(II)}$.
This shows that $\{\alpha_N\}$ satisfies the condition (2.11). Moreover, again by (B.3), 
$$
\Pi_{\theta_N}^N \Bigl(\frac{{\cal U}_N^{(II)}-b_N^{(II)} }{a_N^{(II)}} \geq w_1^{(II)}\Bigr) \leq \Pi_{\theta_N}^N\bigl({\cal U}_N^{(II)} 
\geq u_{\alpha_{N}N}^{(II)}\bigr) \leq
\Pi_{\theta_N}^N \Bigl(\frac{{\cal U}_N^{(II)}-b_N^{(II)} }{a_N^{(II)}} \geq w_0^{(II)}\Bigr).
$$
Taking appropriate limits of both sides we infer that the above chosen sequence $\{\alpha_N\}$, in addition to satisfy (2.11), belongs to $\mathbb{L}^*$, as
$$
0 < 1-U_2^{(II)}(w_1^{(II)}) \leq \liminf_{N \to \infty} \Pi_{\theta_N}^N\bigl({\cal U}_N^{(II)} \geq u_{\alpha_{N}N}^{(II)}\bigr) \leq
\limsup_{N \to \infty} \Pi_{\theta_N}^N\bigl({\cal U}_N^{(II)} \geq u_{\alpha_{N}N}^{(II)}\bigr)$$
$$ \leq 1-U_1^{(II)}(w_0^{(II)}) < 1.
$$ 
\hfill $\Box$


\noindent
{\bf Appendix C: Proof of Theorem 2}\\

\noindent 
The argument follows the idea developed in \cite{r16} and exploits the Koml\'os-Major-Tusn\'ady  inequality for the uniform empirical process. Therefore, consider two probability spaces, two independent sequences $\{B_m^{'}\}$ and $\{B_n^{''}\}$ of Brownian bridges defined on them, and two independent sequences of uniform empirical processes $\{e_m^{'}\}$ and $\{e_n^{''}\}$, defined on the same space, such that for all $m, n$ and $w \in \mathbb{R}$
$$
Pr\Bigl(\sup_{t\in [0,1]} |e_m^{'}(t) - B_m^{'}(t)| \geq \frac{w + C \log m}{\sqrt m}\Bigr) \leq L\exp\{-lw\},
$$
$$
Pr\Bigl(\sup_{t\in [0,1]} |e_n^{''}(t) - B_n^{''}(t)| \geq \frac{w + C \log n}{\sqrt n}\Bigr) \leq L\exp\{-lw\},
$$
where $C, L$ and $l$ are absolute positive constants. On the other hand, 
$$
\sqrt{\frac{mn}{N}} \Bigl\{\hat G_n(z) - \hat F_m(z)\Bigr\} \stackrel{D}{=}\sqrt{\frac{m}{N}} e_n^{''}(J_1(z))-\sqrt{\frac{n}{N}} e_m^{'}(J_1(z)),
$$
where $\stackrel{D}{=}$ denotes the equality in distribution while $B_N^0\stackrel{D}{=} \sqrt{\frac{m}{N}}B_n^{''} - \sqrt{\frac{n}{N}}B_m^{'}$ is a Brownian bridge. Hence, by the above and the property $Pr\bigl(\sup_{t\in[0,1]} B_N^0(t) \geq w\bigr) = \exp\{-2w^2\},\;w \in \mathbb{R},$ we get 
$$
P_0^N({\cal V}_N \geq w_N \sqrt N ) = Pr\Bigl(\sup_{t\in[0,1]} \Big\{\sqrt{\frac{m}{N}} e_n^{''}(t) - \sqrt{\frac{n}{N}} e_m^{'}(t)\Bigr\} \geq w_N \sqrt N \Bigr)
$$
$$
\leq Pr\bigl(\sup_{t\in[0,1]} B_N^0(t) \geq (1-\sqrt{w_N})w_N \sqrt N\bigr) 
$$
$$
+ Pr\Bigl(\sup_{t\in[0,1]} | e_m^{'}(t) - B_m^{'}(t)| \geq \sqrt{\frac{N}{n}}\frac{\sqrt{w_N}}{2} w_N \sqrt N\Bigr)
+ Pr\Bigl(\sup_{t\in [0,1]}| e_n^{''}(t)-B_n^{''}(t)| \geq \sqrt{\frac{N}{m}}\frac{\sqrt{w_N}}{2} w_N \sqrt N\Bigr)
$$
$$
\leq (1+o(1))\exp\bigl\{-2(1-\sqrt{w_N})^2w_N^2 N\bigr\}. 
$$ 
Analogously we obtain $P_0^N({\cal V}_N \geq  w_N \sqrt N) \geq (1+o(1))\exp\bigl\{-2(1+\sqrt{w_N})^2w_N^2 N\bigr\}$.  $\hfill{\Box}$\\

\noindent
{\bf Appendix D: Proof of Theorem 3 and verification of (3.10)}\\

\noindent
Since we like to apply some results of \cite{r15} therefore we have to adjust our statistics to pertaining ones considered in that paper. First of all note that the results of that paper apply, as well, to rank statistics with the score function depending on N.

Next observe that, by (i) and (ii), it holds everywhere
$$
\max_{1 \leq j \leq \Delta(N)} \Bigl|\sum_{i=1}^{N} c_{Ni}\,{\ell}_j\Bigl(\frac{R_i-0.5}{N}\Bigr) - \sum_{i=1}^{N} c_{Ni}\,{\ell}_j\Bigl(\frac{R_i}{N}\Bigr)   \Bigr| \leq
$$
$$
\max_{1 \leq j \leq \Delta(N)} \sqrt{\frac {N}{mn} } \frac{1}{\sqrt{\pi_{jN}(1-\pi_{jN}) }} 
\leq
\frac{\Delta(N) +1}{\sqrt {\Delta(N) }}\sqrt{\frac{N}{mn} } = O\bigl(\sqrt{\frac{\Delta(N)}{N}}\bigr)=o(1).
\eqno(D.1)
$$
So, we can abandon the correction for continuity in ${\cal L}_j$.

Finally, we construct appropriate continuous approximation of the score functions ${\ell}_j$, $j=1,2,..,\Delta(N)$. For this purpose some auxiliary notation are introduced. 

For a fixed $\tau \in (0,1)$ set 
$$
l(t;\tau)=- \sqrt{\frac{1-\tau}{\tau}}\, {\bf 1}(0\leq t <\tau) +
\sqrt{\frac{\tau}{1- \tau}}\, {\bf 1}(\tau \leq t \leq 1).
$$

Given $\epsilon \in (0,\tau(1-\tau))$ we shall modify $l$ on the interval ${\mathbb{I}}_{\epsilon}(\tau) = [\tau (1-\epsilon), \tau (1-\epsilon) + \epsilon]$ containing the jump point $\tau$. To this end introduce the function
$r(t;\tau,\epsilon)$ which is 0 outside ${\mathbb{I}}_{\epsilon}(\tau)$,
$$r(t;\tau,\epsilon)= \sqrt{\frac{1-\tau}{\tau}} + \frac{1}{\epsilon} \sqrt{\frac{1-\tau}{\tau^3}}(t-\tau)\;\;\;\mbox{if}\;\;\;\tau(1-\epsilon)  \leq t < \tau
$$
and
$$
r(t;\tau,\epsilon)= - \sqrt{\frac{\tau}{1-\tau}} + \frac{1}{\epsilon}\sqrt{\frac{\tau}{(1-\tau)^3}}(t-\tau)\,\,\,\mbox{if}\;\;\;\tau \leq t \leq \tau (1-\epsilon) + \epsilon.
$$
Then define 
$$
\bar l (t;\tau,\epsilon)=\sqrt{\frac{3}{3-2\epsilon}}\Bigl[l(t;\tau) + r(t;\tau,\epsilon)\Bigr]
$$
and note that $\bar l (t;\tau,\epsilon)$ is piecewise linear, absolutely continuous, and satisfies  $\int_0^1 \bar l(t;\tau,\epsilon)dt = 0$, and $\int_0^1 \bar l ^2 (t;\tau,\epsilon)dt = 1.$ Moreover, on the interval ${\mathbb{I}}_{\epsilon}(\tau)$ it holds that $|l(t;\tau)- \bar l (t;\tau,\epsilon)| \leq 1/\sqrt{\tau (1-\tau)}$ while outside this interval 
$|l(t;\tau)- \bar l (t;\tau,\epsilon)| \leq \epsilon /\sqrt{\tau (1-\tau)}$.   
For $\epsilon < 1/N$ there is at most one point $R_i/N$ in the interval ${\mathbb{I}}_{\epsilon}(\tau)$. Hence 
$$
\Big|\sum_{i=1}^{N} c_{Ni}\,l\Bigl(\frac{R_i}{N};\tau\Bigr) - \sum_{i=1}^{N} c_{Ni}\,\bar l\Bigl(\frac{R_i}{N};\tau,\epsilon\Bigr)\Bigr| \leq 
\Bigl(\epsilon + \lceil N \epsilon \rceil\Bigr)\frac{1}{\sqrt{\tau (1-\tau)}} \sqrt{\frac{N}{mn}}.
\eqno(D.2)
$$
Take now $\tau =\pi_{jN},\;\epsilon=\epsilon_N = 1/(2N)$ and define 
$$
\bar {\cal L}_j = \bar {\cal L}_{jN}= \sum_{i=1}^{N} c_{Ni}\,\bar l \Bigl(\frac{R_i}{N};\pi_{jN},\epsilon_N \Bigr)\;\;\;\;\;\mbox{and}\;\;\;\;\;
\bar {\cal T}_N = \max_{1 \leq j \leq \Delta(N)}\{-\bar {\cal L}_j\}.
$$
Then, by (D.1) and (D.2), for all $N$ large enough we have everywhere
$$
\bigl|{\cal L}_j - \bar {\cal L}_j\bigr| \leq \frac{2}{\sqrt{\eta (1-\eta)}} \sqrt{{\Delta(N)}/{N}}\;\;\;\;\;\mbox{and}\;\;\;\;\; 
\bigl|{\cal T}_N - \bar {\cal T}_N\bigr| \leq \frac{2}{\sqrt{\eta (1-\eta)}} \sqrt{{\Delta(N)}/{N}}.
\eqno(D.3)
$$

For each of the rank statistic $\bar {\cal L}_j,\;j=1,...,\Delta(N),$ we shall apply Theorem 3.4 of Inglot (2012). Note that in our situation we need to insert there $\Psi(1)$ in place of $\Psi(d(N))$, where $\Psi(1)=\int_0^1\bigl|\frac{\partial}{\partial t} \bar l(t;\pi_{jN},\epsilon_N)\bigr|dt$; cf. (3.13) ibidem. We have
$\Psi(1)=\sqrt{3}/\sqrt{(3-2\epsilon_N)\pi_{jN}(1-\pi_{jN})}$. Moreover, $\lambda_N $ appearing in that theorem equals 1 in our application. The above yields
$$
P_0^N\bigl(\bigl|\bar {\cal L}_j\bigr| \geq w_N \sqrt N \bigr) = 
\exp\Bigl\{-\frac{1}{2}Nw_N^2 + O(Nw_N^{2+\upsilon /2})  + O(\log Nw_N^2)\Bigr\}
\eqno(D.4)
$$
uniformly in $j$. This implies that 
$$
P_0^N\bigl(\bar {\cal T}_N \geq w_N \sqrt N \bigr) = 
\exp\Bigl\{-\frac{1}{2}Nw_N^2 + O(Nw_N^{2+\upsilon /2}) + O(\log Nw_N^2)+ O(\log \Delta(N))\Bigr\}.
\eqno(D.5)
$$
In view of (D.3) and (i), (D.5) yields
$$
P_0^N\bigl( {\cal T}_N \geq w_N \sqrt N \bigr) = 
\exp\Bigl\{-\frac{1}{2}Nw_N^2 + O(Nw_N^{2+\upsilon /2}) + O(\log Nw_N^2)+ O(\log \Delta(N))\Bigr\}.
\eqno(D.6)
$$
Since $\Delta(N)=o(N)$, then by (iii),  $O(\log \Delta(N))+ O(\log N w_N^2)=o(Nw_N^2)$. Hence (3.12) follows. $\hfill{\Box}$\\

\noindent{\bf Verification of (3.10)}\\

\noindent
We argue similarly as in the proof of Lemma A.1 in \cite{r25}.

Put 
$$ 
\widetilde{\cal L}_j=-\sqrt{N/mn}\,{\cal L}_j.
$$ 
Let $Z_1,...,Z_N$ denote the pooled sample $X_1,...,X_m,Y_1,...,Y_n$ and let $Z_{(r)}$ stand for the $r$-th order statistic of the pooled sample. 
For any $j=1,...,\Delta(N)$ we have
$$\widetilde{\cal L}_j=-\int_{Z_{(1)}}^{\infty} l_j\left(\hat{J}_N(x)-\frac{1}{2N}\right)\,d(\hat{G}_n(x)-\hat{F}_m(x))
=-\int_{1/N}^1l_j\left(t-\frac{1}{2N}\right)\,d(\hat{G}_n-\hat{F}_m)\circ\hat{J}_N^{-1}(t).$$
Applying to the last expression the integration by parts formula, cf. (1) in \cite{r101}, p. 115, we get 
$$
\widetilde{\cal L}_j=
(\hat{G}_n-\hat{F}_m)\circ\hat{J}_N^{-1}\left(\frac{1}{N}\right)\,l_j\left(\frac{1}{2N}\right)+
\int_{1/N}^1 (\hat{G}_n-\hat{F}_m)\circ\hat{J}_N^{-1}(t)\,d\,l_j\left(t-\frac{1}{2N}\right)
$$
$$
=(\hat{G}_n-\hat{F}_m)\circ\hat{J}_N^{-1}\left(\frac{1}{N}\right)\,l_j\left(\frac{1}{2N}\right)+
\frac{1}{\sqrt{\pi_{jN}(1-\pi_{jN})}}(\hat{G}_n-\hat{F}_m)\circ\hat{J}_N^{-1}\left(\pi_{jN}+\frac{1}{2N}\right).
$$
Set 
$$
\widetilde{W}_j=\frac{1}{\sqrt{\pi_{jN}(1-\pi_{jN})}}(\widehat{G}_n-\widehat{F}_m)\circ \widehat{J}_N^{-1}(\pi_{jN}).
$$ 
Then we have
$$
|\widetilde{\cal L}_j-\widetilde{\cal W}_j|\leq\frac{1}{\sqrt{\pi_{jN}(1-\pi_{jN})}}\left[\left|(\hat{G}_n-\hat{F}_m)\circ\hat{J}_N^{-1}\left(\frac{1}{N}\right)
\max\{\pi_{jN},1-\pi_{jN}\}\right|\right.
$$
$$
\left.+\left|(\hat{G}_n-\hat{F}_m)\circ\hat{J}_N^{-1}\left(\pi_{jN}+\frac{1}{2N}\right)-(\hat{G}_n-\hat{F}_m)\circ\hat{J}_N^{-1}(\pi_{jN})\right|\right].
\eqno (D.7)
$$

By the definition, $\widehat{J}_N^{-1}(t)$ equals $Z_{(\lceil N\,t\rceil)},\;t\in (0,1)$. So, the first term in (D.7) is majorized by 
$$ \frac{1}{\sqrt{\pi_{jN}(1-\pi_{jN})}} \max\left\{\frac{1}{m},\frac{1}{n}\right\}. \eqno (D.8)$$
When $\lceil N\pi_{jN}+0.5\rceil=\lceil N\pi_{jN}\rceil$ then the second term in (D.7) equals 0.
When $\lceil N\pi_{jN}+0.5\rceil=\lceil N\pi_{jN}\rceil+1$ then the second term in (D.7) is also majorized by (D.8)
Hence
$$
|\widetilde{\cal L}_j-\widetilde{\cal W}_j|\leq\frac{2}{\sqrt{\pi_{jN}(1-\pi_{jN})}}\max\left\{\frac{1}{m},\frac{1}{n}\right\}.
$$
Since 
$$
{\cal T}_N=\sqrt{\frac{mn}{N}}\max_{1\leq j\leq \Delta (N)}\widetilde{L}_j\;\;\;\mbox{and}\;\;\; {\cal W}_N=\sqrt{\frac{mn}{N}}\max_{1\leq j\leq \Delta (N)}\widetilde{W}_j, 
$$
by the triangle inequality, we get
$$
|{\cal T}_N-{\cal W}_N|\leq \sqrt{\frac{mn}{N}}\max_{1\leq j\leq \Delta (N)}|\widetilde{\cal L}_j-\widetilde{\cal W}_j|
$$
$$
\leq \frac{1}{\sqrt{N}}\max\left\{\sqrt{\frac{1-\eta_N}{\eta_N}},\sqrt{\frac{\eta_N}{1-\eta_N}}\right\} \times
\frac{2}{\min\{\sqrt{\pi_{1N}(1-\pi_{1N})},\sqrt{\pi_{\Delta(N)N}(1-\pi_{\Delta(N)N})}\}}.
$$
This, after elementary argument, yields (3.10). $\hfill{\Box}$\\

\noindent
{\bf Appendix E: Proof of Theorem 4}\\ 

\noindent
We have $\sqrt m \{\hat F_m - F_{1N}\} \stackrel{D}{=} e_m^{(1)} \circ F_{1N},\;\;\sqrt n \{\hat G_n - G_{1N}\} \stackrel{D}{=} e_n^{(2)} \circ G_{1N}$, where $e_m^{(1)}$ and $e_n^{(2)}$ are independent uniform empirical processes defined on an appropriate probability space. In particular, one can use the KMT constructions applied  in the proof of Theorem 2.  Hence
$$
\sqrt{\frac{mn}{N}} \bigl\{\hat G_n - \hat F_m\bigr\} \stackrel{D}{=} \sqrt{\frac{m}{N}} e_n^{(2)} \circ G_{1N} - \sqrt{\frac{n}{N}} e_m^{(1)} \circ F_{1N} + \sqrt{\frac{mn}{N}} \vartheta_N (G_1 - F_1)
\eqno(E.1)
$$ 
and
$$
\Pi_{\vartheta_N}^N({\cal V}_N - b_{{\cal V}}(\Pi_{\vartheta_N}^N) \leq w) \leq 
Pr \Bigl(\sqrt{\frac{m}{N}} e_n^{(2)} \circ G_{1N}(z_0) - \sqrt{\frac{n}{N}} e_m^{(1)} \circ F_{1N}(z_0) \leq w\Bigr).
\eqno(E.2)
$$
Since $G_{1N}(z_0) \to J_1(z_0)$ and $F_{1N}(z_0) \to J_1(z_0)$, therefore the random variable on the right hand side of (E.2) has asymptotic $N(0, \sqrt{J_1(z_0)[1-J_1(z_0)}])$ law.  
This justifies the form of $V_2$.

On the other hand, by (E.1) we infer that 
$$
\Pi_{\vartheta_N}^N \Bigl({\cal V}_N - b_{\cal V} (\Pi_{\vartheta_N}^N) \leq w\Bigr) \geq
\Pi_{\vartheta_N}^N \Bigl(\sqrt{\frac{mn}{N}} \sup_{z \in \mathbb{R}}\Bigl\{\hat G_n(z) - \hat F_m(z) - \vartheta_N [G_1(z)-F_1(z)]\Bigr\} \leq w\Bigr) =
$$
$$
Pr\Bigl(\sup_{z \in \mathbb{R}} \Bigl\{ \sqrt{\frac{m}{N}} e_n^{(2)} \circ G_{1N}(z) - \sqrt{\frac{n}{N}} e_m^{(1)} \circ F_{1N}(z)\Bigr\} \leq w\Bigr).
$$
Since $e_m^{(1)} \circ F_{1N} \Rightarrow B^{(1)}$ and $e_n^{(2)} \circ G_{1N} \Rightarrow B^{(2)}$, where $B^{(1)}$ and $B^{(2)}$ are independent Brownian bridges while $\Rightarrow$ denotes weak convergence, the form of $V_1$ follows. $\hfill{\Box}$\\

\noindent
{\bf Appendix F: Proof of Theorem 5}\\

\noindent
{\bf F.1. Preliminaries.}
We shall prove (3.16) and (3.17) for the statistic ${\cal W}_N$. Since (ii) implies that $\Delta(N) = o(\sqrt N / \log N)$, therefore (3.10) justifies such approach.

Set 
$$\kappa_N =\sqrt{\eta_N(1-\eta_N)}\;\;\;\;\mbox{and}\;\;\;\;\nu_N=\eta_N - \eta $$ 
and introduce two auxiliary processes on $\mathbb{R}$
$$
\zeta_N(z)=\sqrt N \kappa_N[\hat G_n(z)-G_{1N}(z)]-\sqrt N \kappa_N[\hat F_m(z)-F_{1N}(z)],
$$
$$
\xi_N(z)=\sqrt N (1-\eta_N)[\hat G_n(z)-G_{1N}(z)]+\sqrt N \eta_N[\hat F_m(z)-F_{1N}(z)].
$$
For $\hat J_N(z)=\eta_N \hat F_m(z)+(1-\eta_N) \hat G_n(z)$ put 
$$\hat z_{jN}=\hat J_N^{-1}(\pi_{jN}).$$ 
Additionally, set $V_N(z)={\sqrt N}\kappa_N\{\hat G_n(z)-\hat F_m (z)\}$. With these notation 
$$ 
\hat J_N(z)=\frac{1}{\sqrt N} \xi_N(z) +J_1(z) -\theta_N \nu_N \bar A(J_1(z))
\eqno(F.1)
$$ 
while  
$$
V_N(z)=\zeta_N(z) + \sqrt N \theta_N \kappa_N \bar A(J_1(z))\;\;\;\mbox{and}\;\;\;
{\cal W}_N =\max_{1 \leq j \leq \Delta(N)}\frac{V_N(\hat z_{jN})}{\sqrt{\pi_{jN}(1-\pi_{jN})}}.
\eqno(F.2)
$$
Now, let us reperametrize $F_{1N}$ and $G_{1N}$ in ({3.2}) to a classical form in the two-sample scheme, which we shall exploit below. For $t\in (0,1)$ set
$$
\bar A(t) = \bar A(t;\eta) = (G_1-F_1) \circ J_1^{-1}(t),\ \ \bar A^+ (t) = \bar A^+ (t;\eta) = \max\{\bar A(t),0\}.
\eqno (F.3)$$
With the above notation,  (3.2) can be written as
$$
F_{1N}= J_1 - \vartheta_N (1-\eta) \bar A \circ J_1,\;\;G_{1N}=J_1 + \vartheta_N \eta \bar A \circ J_1.
\eqno (F.4)$$
By ({F.3}) it follows that $\bar A(t)$ is absolutely continuous and $F_1(z) = J_1(z) - (1-\eta)\bar A \circ J_1(z),\;\;G_1(z)=J_1(z)+\eta \bar A \circ J_1(z)$. Hence,  $\bar a(t)=\bar A'(t)$ exists almost everywhere (with respect to the Lebesgue measure) and it holds that 
$$
-\eta ^{-1} \leq \bar a(t) \leq (1-\eta)^{-1} 
\;\;\;\mbox{and}\;\;\int_0^1\bar a(t)dt =0;
$$
cf. Behnen and Neuhaus (1983,1989), for example. Note also that 
$$
\frac{d F_{1N}}{d J_1} \circ J_1^{-1}(t) = 1 -\vartheta_N (1-\eta)\bar a(t)\;\;\mbox{and}\;\;\frac{d G_{1N}}{d J_1} \circ J_1^{-1}(t) = 1 +\vartheta_N \eta \bar  a(t).
\eqno (F.5)$$
In consequence, for each $\eta \in (0,1)$ we have $\bar A(0;\eta)=\bar A(1;\eta)=0$ and
$$
\lim_{t \to 0+}A^*(t;\eta)=\lim_{t \to 1-}A^*(t;\eta)=0,\;\;\;\;\mbox{where}\;\;A^*(t;\eta)=\frac{\bar A(t;\eta)}{\sqrt{t(1-t)}}.
\eqno (F.6)$$
By ({F.6}) there exists $\delta \in (0,1/2)$ such that 
$$
\max_{j : \pi_{jN} \notin [2\delta,1-2\delta]} \frac{\bar A(\pi_{jN})}{\sqrt{\pi_{jN}(1-\pi_{jN})}} \leq
\frac{1}{2} \mu_0,\;\;\;\mbox{where}\;\;\;\mu_0=
\max_{j : \pi_{jN} \in [2\delta,1-2\delta]} \frac{\bar A(\pi_{jN})}{\sqrt{\pi_{jN}(1-\pi_{jN})}}.
\eqno(F.7)
$$
To increase readability of the proof of (3.16) and (3.17) we formulate now some partial results, which we shall justify at Subsections F.2 - F.6.

For $\delta$ defined via (F.7) set 
$$
z_1=J_1^{-1}(\delta)\;\;\;\mbox{and}\;\;\;z_2=J_1^{-1}(1-\delta).
$$ 
Recall that $\Pi_{\theta_N}^N = P_{\theta_N}^{m(N)} \times Q_{\theta_N}^{n(N)}$, where $P_{\theta_N}$ and $Q_{\theta_N}$ are defined via (F.4) with $\theta_N \to 0$ and $N \theta_N^2 \to \infty$. In the succeeding lemmas we specify sufficient conditions on $\{\theta_N\}$ for them to hold. Throughout $C$ is an absolute constant, not necessarily the same in all places.\\

\noindent{\bf Lemma F.1.}\\
(a) {\it If $\theta_N \to 0$ and (ii) holds then
$$
\lim_{N \to \infty} \Pi_{\theta_N}^N\Bigl(\sup_{z \in [z_1,z_2]} \frac{|\zeta_N(z)|}{\sqrt {\hat J_N(z)[1-\hat J_N(z)]}} \leq w\Bigr) 
= Pr \Bigl(\sup_{t \in [\delta,1-\delta]} \frac{|B(t)|}{\sqrt{t(1-t)}} \leq w\Bigr),\;w \in \mathbb{R}_{+},
\eqno(F.8)
$$
where $B$ is a Brownian bridge.}
\\
(b) {\it Assume that (i), (ii) and (iii) of Theorem 5 hold.
Set 
$$
\mathbb{E}_{0N} = \Bigl\{\sup_{z \in \mathbb{R}}|\hat J_N(z) - J_1(z)| \leq \sqrt{\frac{\log N}{N}}\Bigr\},\;\;\;
\mathbb{E}_{1N} = \Bigl\{\sup_{z \in \mathbb{R}}|\hat J_N(z) - J_1(z)| \leq \frac{C}{\theta_N \Delta(N) \sqrt N}\Bigr\},
$$
and
$$
\mathbb{E}_{2N} =\Bigl\{ \max_{1 \leq j \leq \Delta(N)}\frac{|\bar A(\hat J_N(\hat z_{jN}))- \bar A(J_1(\hat z_{jN}))|}{\sqrt{\pi_{jN}(1-\pi_{jN})} }  \leq \frac{1}{\theta_N \sqrt{N \Delta(N)}}  \Bigr\}.
$$
Then }
$$
\lim_{N \to \infty} \Pi_{\theta_N}^N \Bigl(\mathbb{E}_{0N}\Bigr)=
\lim_{N \to \infty} \Pi_{\theta_N}^N \Bigl(\mathbb{E}_{1N}\Bigr)= \lim_{N \to \infty} \Pi_{\theta_N}^N \Bigl(\mathbb{E}_{2N}\Bigr)=1.
\eqno(F.9)
$$

\noindent 
Moreover, the following useful bounds take place.
On $\mathbb{E}_{0N}$ we have
$$
|J_1(\hat z_{jN}) - \pi_{jN}| \leq |J_1(\hat z_{jN}) -J_N(\hat z_{jN})| + |J_N(\hat z_{jN})-\pi_{jN}| \leq \sqrt{\frac{\log N}{N}} + \frac{1}{N} 
\leq 2\sqrt{\frac{\log N}{N}}
\eqno(F.10)
$$
while on $\mathbb{E}_{1N}$ 
$$
|J_1(\hat z_{jN}) - \pi_{jN}| \leq 
\frac{C}{\theta_N \Delta(N) \sqrt N} .
\eqno(F.11)
$$

Further introduce 
$$
l_N=J_1^{-1}\Bigl(\frac{\log N}{\sqrt N}\Bigr)\;\;\;\mbox{and}\;\;\;u_N=J_1^{-1}\Bigl(1-\frac{\log N}{\sqrt N}\Bigr)
$$
and note that for $N$ large enough it holds $l_N \leq u_N$.
\\

\noindent
{\bf Lemma F.2.} {\it Suppose that $N\theta_N^2 /\log^2 N  \to \infty$, $\theta_N \sqrt N \nu_N = O(1)$, and $w_N \asymp \theta_N \sqrt N$. Then }
$$
\lim_{N \to \infty} \Pi_{\theta_N}^N\Bigl( \sup_{z \in [l_N,u_N]} \frac{|\zeta_N(z)|}{\sqrt {\hat J_N(z)[1-\hat J_N(z)]}} \leq w_N\Bigr) = 1. 
\eqno(F.12)
$$
\\

\noindent
{\bf Lemma F.3}. {\it Under (i), (ii) and (iii) of Theorem 5, for $\mathbb{E}_{3N}$ given by
$$
\mathbb{E}_{3N} =\Bigl\{{\cal W}_N \leq \max_{j : \hat z_{jN} \in [z_1,z_2]} \frac{V_N(\hat z_{jN})}{\sqrt{\pi_{jN}(1-\pi_{jN})}}\Bigr\}
$$
it holds that }
$$
\lim_{N \to \infty} \Pi_{\theta_N}^N \Bigl(\mathbb{E}_{3N}\Bigr)=
\lim_{N \to \infty} \Pi_{\theta_N}^N \Bigl(\max_{j : \hat z_{jN} \notin [z_1,z_2]} \frac{V_N(\hat z_{jN})}{\sqrt{\pi_{jN}(1-\pi_{jN})}} \leq
\max_{j : \hat z_{jN} \in [z_1,z_2]} \frac{V_N(\hat z_{jN})}{\sqrt{\pi_{jN}(1-\pi_{jN})}}\Bigr) = 1.
\eqno(F.13)
$$

By the above, to prove (3.16) and (3.17) it is enough to consider
$$
\Pi_{\theta_N}^N\Bigl(\Bigl\{{\cal W}_N - b_{{\cal T}}(\Pi_{\theta_N}^N) \leq w\Bigr\}\cap \bigcap _{j=0}^3 \mathbb{E}_{jN}\Bigr).
$$
\\

\noindent
{\bf F.2. Proof of (3.16).}
Let $j_0 = j_0(N)$ be any index $j$ such that 
$$
\max_{1 \leq j \leq \Delta(N)} \frac{\bar A(\pi_{jN})}{\sqrt{\pi_{jN}(1-\pi_{jN})}} =
\frac{\bar A(\pi_{j_0 N})}{\sqrt{\pi_{j_0 N}(1-\pi_{j_0 N})}}
$$
By (F.7), without loss of generality we can assume that $j_0$ is such that for each $N$ it holds that $\pi_{j_0 N} \in [2\delta,1-2\delta]$.
With this notation 
$$
b_{\cal T}(\Pi_{\theta_N}^N) = \sqrt N \theta_N \kappa_N \frac{\bar A(\pi_{j_0 N})}{\sqrt{\pi_{j_0 N}(1-\pi_{j_0 N})}}.
$$
By (F.6), (F.11) and (i), on the set $\mathbb{E}_{1N}$ 
$$
{\cal W}_N - b_{\cal T}(\Pi_{\theta_N}^N) \geq \frac{\zeta_N(\hat z_{j_0 N}) + \sqrt N \theta_N \kappa_N [\bar A(J_1(\hat z_{j_0 N})) - \bar A(\pi_{j_0 N})] }
{\sqrt{\pi_{j_0 N}(1-\pi_{j_0 N})}}=
\frac{\zeta_N(\hat z_{j_0 N})}{\sqrt{\pi_{j_0 N}(1-\pi_{j_0 N})}} + o(1).
$$

Therefore, to conclude the proof of (3.16) it is enough to show that
$$
\frac{\zeta_N(\hat z_{j_0 N})}{\sqrt{\pi_{j_0 N}(1-\pi_{j_0 N})}}
\stackrel{D}{\longrightarrow} N(0,1)\;\;\;\mbox{as}\;\;\;N \to \infty.
\eqno(F.14)
$$
The main difficulty in proving (F.14) lies in that $j_0=j_0(N)$ may be not unique and changes with $N$. Therefore, we proceed as follows.
When $N$ is growing then, by (i), the partition is getting more dense. Hence, the set of accumulation points of the sequence $\{\pi_{j_0 N}\}$ is nonempty and is contained in $[2\delta,1-2\delta]$. Therefore, it is enough to prove (F.14) for any concentration point and pertaining subsequence of 
$\{\pi_{j_0 N}\}$ converging to it. Set $t_0$ to be any concentration point of the sequence and denote by 
$\{\pi_{j'_0 N'}\}$, $j'_0 =j_0(N')$, a subsequence converging to $t_0$. 
By the definition of $\hat z_{jN}$ and (F.10), on $\mathbb{E}_{0N}$ it holds 
$$
|J_1(\hat z_{j'_0 N'}) -t_0| \leq |J_1(\hat z_{j'_0 N'})- \hat J_{N'}(\hat z_{j'_0 N'})| + |\hat J_{N'}(\hat z_{j'_0 N'}) - \pi_{j'_0 N'}| + |\pi_{j'_0 N'} - t_0|
\leq 2\sqrt{\frac{\log N}{N}} + |\pi_{{j'_0}N'} - t_0|.
$$
and yields 
$$
J_1(\hat z_{j'_0 N'}) \stackrel{\Pi_{\theta_{N'}}^{N'}}{\longrightarrow} t_0.
$$
This and the continuity of $J_1$, imply that subsequence  $\{\hat z_{{j'_0} N'}\}$  converges in $\Pi_{\theta_{N'}}^{N'}$ to $J_1^{-1}(t_0).$ Hence, 
weak convergence of the process $\zeta_N (z)/\sqrt{\hat J_N (z)[1-\hat J_N (z)]}$, to the process\\ $B(J_1(z))/\sqrt{J_1(z)[1-J_1(z)]},\;z \in [z_1,z_2],$
cf. the proof of Lemma F.1,    implies that, under $\Pi_{\theta_{N'}}^{N'}$,
$$
\frac{\zeta_{N'}(\hat z_{j'_0 N'})}{\sqrt{\pi_{j'_0 N'}(1-\pi_{j'_0 N'})}}
\stackrel{D}{\longrightarrow} \frac{B(t_0)}{\sqrt{t_0(1-t_0)}}\;\;\;\mbox{as}\;\;\;N' \to \infty.
\eqno(F.15)
$$
\\
This shows that for any convergent subsequence $\{\pi_{j'_0 N'}\}$ of the sequence $\{\pi_{j_0 N}\}$ the sequence of random variables in (F.15) converges to the same limiting N(0,1) law. This proves (3.16). \hfill${\Box}$\\

\noindent
{\bf F.3. Proof of (3.17).}
Recall that we can restrict attention to $\bigcap_{j=0}^3 \mathbb{E}_{jN}$. In particular, on $\mathbb{E}_{1N} \cap \mathbb{E}_{3N}$, by (F.6), (F.11) and Lipschitz condition for $\bar A$,  we have
$$
{\cal W}_N - b_{\cal T}(\Pi_{\theta_N}^N) \leq 
\max_{j : \hat z_{jN} \in [z_1,z_2]} \frac{V_N (\hat z_{jN})-\sqrt N \theta_N \kappa_N \bar A(\pi_{jN})}{\sqrt{\pi_{jN}(1-\pi_{jN})}}=
\max_{j:\hat z_{jN} \in [z_1,z_2]} \Bigl\{\frac{\zeta_N(\hat z_{jN})}{\sqrt{\pi_{jN}(1-\pi_{jN})}} +$$
$$ 
\frac{\sqrt N \theta_N \kappa_N [\bar A (J_1(\hat z_{jN})) -\bar A (\pi_{jN})]}{\sqrt{\pi_{jN}(1-\pi_{jN})}}\Bigr\} 
\leq\max_{j:\hat z_{jN} \in [z_1,z_2]} \Bigl\{\frac{|\zeta_N(\hat z_{jN})|}{\sqrt{\pi_{jN}(1-\pi_{jN})}}\Bigr\} +  \frac{C}{\sqrt{\Delta(N)}}.
\eqno(F.16)
$$

Now observe that the property $|\hat J_N(\hat z_{jN}) - \pi_{jN}| < 1/N$ implies that 
$|\hat J_N(\hat z_{jN})[1-\hat J_N(\hat z_{jN})] - \pi_{jN}(1-\pi_{jN})| \leq 1/N.$ Hence, by (i), for $N$ large enough 
$$
1-2\frac{\Delta(N)}{N} \leq \sqrt{\frac{ \hat J_N(\hat z_{jN})[1-\hat J_N(\hat z_{jN})]}{\pi_{jN}(1-\pi_{jN})}} \leq 1+2\frac{\Delta(N)}{N}.
\eqno(F.17)
$$
Hence, the right hand side of (F.16) is majorized by 
$$
[1+o(1)] \sup_{z \in [z_1,z_2]}\frac{|\zeta_N(z)|}{\sqrt{\hat J_N(z)[1-\hat J_N(z)]}} + o(1).
$$
By (F.8) of Lemma F.1 the proof is concluded. $\hfill{\Box}$
\\

\noindent
{\bf F.4. Proof of Lemma F.1.}

(a) As in the proof of Theorem 4, an application of strong approximation technique implies that $\zeta_N \Rightarrow B^{(1)}\circ J_1$  and $\xi_N \Rightarrow B^{(2)}\circ J_1$, where $B^{(1)}$ and $B^{(2)}$ are independent Brownian bridges. Moreover, (F.1) implies that $\hat J_N \stackrel{\Pi_{\theta_N}^N}{\longrightarrow} J_1$. Hence (F.8) follows. 

(b) By (F.1) it holds
$$
\Pi_{\theta_N}^N(\mathbb{E}_{0N}^c) \leq \Pi_{\theta_N}^N \Bigl(\sup_{z \in \mathbb{R}} |\xi_N(z) -\theta_N \nu_N \sqrt N  \bar A(J_1(z))| > \sqrt{\log N}\Bigr) 
$$
$$
\leq \Pi_{\theta_N}^N\Bigl(\sup_{z \in \mathbb{R}}|\xi_N(z)| \geq \sqrt {\log N} - C  \theta_N \nu_N \sqrt N \Bigr)
$$
and, by the weak convergence of $\xi_N$ and the assumption (iii), $\lim_{N \to \infty} \Pi_{\theta_N}^N(\mathbb{E}_{0N}^c)=0$.\\
Analogously,
$$
\Pi_{\theta_N}^N(\mathbb{E}_{1N}^c) \leq \Pi_{\theta_N}^N \Bigl(\sup_{z \in \mathbb{R}} |\xi_N(z) -\theta_N \sqrt N \nu_N \bar A(J_1(z))| > \frac{C}{\theta_N \Delta (N)}\Bigr)  
$$
$$
\leq \Pi_{\theta_N}^N\Bigl(\sup_{z \in \mathbb{R}}|\xi_N(z)| \geq \frac{C}{\theta_N \Delta (N)} + O(1)\Bigr)
$$
and we infer that $\Pi_{\theta_N}^N(\mathbb{E}_{1N}^c) \to 0$. Moreover, since $\bar A$ is Lipschitz one, then, with some $C$, we have 
$$
\Pi_{\theta_N}^N(\mathbb{E}_{2N}^c) \leq \Pi_{\theta_N}^N \Big(\frac{\sup_{z \in \mathbb{R}}|\hat J_N(z) - J_1(z)|}{\min_{1 \leq j \leq \Delta(N)} \sqrt{\pi_{jN } (1-\pi_{jN})}} >
\frac{C }{\theta_N \sqrt{N \Delta(N)} }\Bigr).
$$
Thus, by (i), it holds $\Pi_{\theta_N}^N(\mathbb{E}_{2N}^c) \leq \Pi_{\theta_N}^N(\mathbb{E}_{1N}^c)$, with some appropriate $C$ in $\mathbb{E}_{1N}$. Hence, the proof of (F.9) is completed. $\hfill{\Box}$\\

\noindent 
{\bf F.5. Proof of Lemma F.2.}
On $\mathbb{E}_{0N}$, given in Lemma B.1, for $z \in (0,1)$ it holds
$|\hat J_N(z)[1-\hat J_N(z)] -J_1(z)[1-J_1(z)]| \leq \sqrt{\log N /N}$. Hence, by the definition of $l_N$ and $u_N$, for $z \in [l_N,u_N]$
$$
\Bigl|\frac{\hat J_N(z)[1-\hat J_N(z)]}{J_1(z)[1-J_1(z)]}-1 \Bigr|\leq \frac{4}{\sqrt{\log N}}.
$$
For $N$ large enough this implies 
$$
\Pi_{\theta_N}^N\Bigl( \sup_{z \in [l_N,u_N]} \frac{|\zeta_N(z)|}{\sqrt {\hat J_N(z)[1-\hat J_N(z)]}} \geq w_N\Bigr)
$$
$$
\leq \Pi_{\theta_N}^N\Bigl( \sup_{z \in [l_N,u_N]} \frac{|\zeta_N(z)|}{\sqrt { J_1(z)[1-J_1(z)]}} \geq w_N \bigl(1-\frac{4}{\sqrt{\log N}}\bigr)\Bigr) + o(1).
\eqno(F.18)
$$
As in the proof of Theorem 2, consider now the uniform empirical process $e_{m}^{'}$ and related Brownian Bridge $B_{m}^{'}$ and independent on them $e_{n}^{''}$ and $B_{n}^{''}$ such that the KMT inequalities hold for them. Under $\Pi_{\theta_N}^N$, $\sqrt m [\hat F_m - F_{1N}] \stackrel{D}{=}e_m^{'}(F_{1N})$,   
 $\sqrt n [\hat G_n - G_{1N}] \stackrel{D}{=}e_n^{''}(G_{1N})$ and $\zeta_N  \stackrel{D}{=} \tilde \zeta_N =\sqrt{\frac{m}{N}} e_n^{''}(G_{1N}) - \sqrt{\frac{n}{N}} e_m^{'}(F_{1N})$. 
Set $B_N=\sqrt{\frac{m}{N}}B_m^{''} -\sqrt{\frac{n}{N}}B_n^{'}$. Then $B_N$ is a Brownian bridge. Therefore,
we can majorize the first component of (F.18) as follows
$$
Pr \Bigl( \sup_{z \in [l_N,u_N]} \frac{|B_N (J_1(z))|}{\sqrt { J_1(z)[1-J_1(z)]}} \geq w_N \bigl(\frac{1}{2}-\frac{4}{\sqrt{\log N}}\bigr)\Bigr)
$$
$$ +
Pr \Bigl( \sup_{z \in [l_N,u_N]} \frac{|\tilde \zeta_N(z) - B_N(J_1(z))|}{\sqrt { J_1(z)[1-J_1(z)]}} \geq \frac{w_N}{2}\Bigr).
\eqno(F.19)
$$
Due to the definition of $l_N, u_N$, the assumptions $w_N \asymp \theta_N \sqrt N, \theta_N \sqrt N / \log N \to \infty$, Darling and Erd\H{o}s result, cf. Lemma 4.4.1 in \cite{r100}, implies that the first component of (F.19) tends to 0.

Using again the form of $l_N$ and $u_N$, the second component of (F.19) for large $N$ is majorized by 
$$
Pr\bigl(\sup_{z \in \mathbb{R}} |\tilde \zeta_N(z) - B_N(J_1(z))| \geq w_N^*\bigr),
\eqno(F.20)
$$
where $w_N^*=w_N \sqrt{\log N}/(4N^{1/4})$. The structure of $\tilde \zeta_N$, given above, allows to majorize (F.20) as follows
$$
Pr\bigl(\sup_{z \in \mathbb{R}}|e_m^{'}(F_{1N}(z))-B_m^{'}(J_1(z))| \geq \kappa_N w_N^*\bigr) +
Pr\bigl(\sup_{z \in \mathbb{R}}|e_n^{''}(G_{1N}(z))-B_n^{''}(J_1(z))| \geq \kappa_N w_N^*\bigr) \leq
$$
$$
Pr\bigl(\sup_{t \in (0,1)}|e_m^{'}(t)-B_m^{'}(t)| \geq \kappa_N w_N^*/2\bigr) +
Pr\bigl(\sup_{t \in (0,1)}|e_n^{''}(t)-B_n^{''}(t)| \geq \kappa_N w_N^*/2\bigr) +
\eqno(F.21)
$$
$$
2Pr\bigl(\sup_{0 \leq t \leq 1-C\theta_N} \sup_{0 \leq h \leq C\theta_N} |B_m^{'}(t+h)-B_m^{'}(t)| \geq \kappa_N w_N^*/2\bigr).
$$
Since $\eta_N \to \eta, w_N \asymp \theta_N \sqrt N$ and $\theta_N \sqrt N \to \infty$ an application of the KMT inequality to the two first components of (F.21) shows that these terms are negligible. The last term of (F.21) requires standard analysis of increments of the Brownian bridge. Applying for this purpose Lemma A of \cite{r15} with $h=C\theta_N, y=y_N=\kappa_N w_N^*/(2\sqrt{C\theta_N})$ and $\delta=1/4$ finishes the proof, as $y \to \infty$ faster than $\log N$. 
$\hfill{\Box}$\\

\noindent
{\bf F.6. Proof of Lemma F.3.}
Recall that $v(t)=\sqrt{t(1-t)}$. Note that
$$
\mathbb{E}_{3N}^c \subset \Bigl\{\max_{j: \hat z_{jN} \notin [z_1,z_2]}\frac{V_N (\hat z_{jN})}{v(\pi_{jN})} > 
 \max_{j: \hat z_{jN} \in [z_1,z_2]}\frac{V_N (\hat z_{jN})}{v(\pi_{jN})}\Bigr\}.
$$

Throughout we restrict attention to $\mathbb{E}_{0N}$. By (F.10), for sufficiently large $N$,\\ $\{j: z_1 \leq \hat z_{jN} \leq z_2\} = \{j: \delta \leq J_1(\hat z_{jN}) \leq 1-\delta\} \supset \{j: 2\delta \leq \pi_{jN} \leq 1-2\delta\}$. Hence $\{j: \hat z_{jN} \notin [z_1,z_1]\} \subset \{j: \pi_{jN} \notin [2\delta, 1-2\delta]\}$ and 
$$
\max_{j: \hat z_{jN} \notin [z_1,z_2]}\frac{V_N (\hat z_{jN})}{v(\pi_{jN})} \leq 
\max_{j: \hat z_{jN} \notin [z_1,z_2]}\frac{V_N (\hat z_{jN}) - \sqrt N \theta_N \kappa_N \bar A(\pi_{jN})}{v(\pi_{jN})} +
\max_{j: \pi_{jN} \notin [2\delta,1-2\delta]}\frac{\sqrt N \theta_N \kappa_N \bar A (\pi_{jN})}{v(\pi_{jN})}.
\eqno(F.22)
$$
Using (F.6), (F.7), (F.10) and (iii) we conclude
$$
\max_{j: \hat z_{jN} \notin [z_1,z_2]}\frac{V_N (\hat z_{jN})}{v(\pi_{jN})} \leq 
\max_{j: \hat z_{jN} \notin [z_1,z_2]}\frac{|V_N (\hat z_{jN}) - \sqrt N \theta_N \kappa_N \bar A(\pi_{jN})|}{v(\pi_{jN})} +
\frac{1}{2}\sqrt N \theta_N \kappa_N \mu_0 
$$
$$
\leq \max_{j: \hat z_{jN} \notin [z_1,z_2]}\frac{|\zeta_N(\hat z_{jN})|}{v(\pi_{jN})} +\rho_N^{(2)},
\eqno(F.23)
$$
where 
$\rho_N^{(2)}=\sqrt N \theta_N \kappa_N (\mu_0/2 + C \sqrt{\Delta(N) \log N /N}) \asymp \sqrt N \theta_N.$

Analogously,
$$
\max_{j: \hat z_{jN} \in [z_1,z_2]}\frac{V_N (\hat z_{jN})}{v(\pi_{jN})} \geq 
\max_{j: \pi_{jN} \in [2\delta, 1-2\delta]}\frac{\sqrt N \theta_N \kappa_N \bar A (\pi_{jN})}{v(\pi_{jN})}-
\max_{j: \hat z_{jN} \in [z_1,z_2]}\frac{|V_N (\hat z_{jN}) - \sqrt N \theta_N \kappa_N \bar A(\pi_{jN})|}{v(\pi_{jN})} 
$$
$$
\geq \rho_N^{(1)} - \max_{j: \hat z_{jN} \in [z_1,z_2]}\frac{|\zeta_N(\hat z_{jN})|}{v(\pi_{jN})},
\eqno(F.24)
$$
where $\rho_N^{(1)}=\sqrt N \theta_N \kappa_N (\mu_0 - C \sqrt{\Delta(N) \log N /N}) \asymp \sqrt N \theta_N.$

The above implies that 
$$
\mathbb{E}_{0N} \cap \mathbb{E}_{3N}^c \subset \Bigl\{\max_{ 1 \leq j \leq \Delta(N)}\frac{|\zeta_N(\hat z_{jN})|}{\sqrt{\pi_{jN}(1-\pi_{jN})}}
> \frac{1}{2}\bigl[\rho_N^{(1)} - \rho_N^{(2)}\bigr] \Bigr\}.
$$
Now, observe that, by  (i), (ii) and (iii), it follows that
$$
\mathbb{E}_{0N} \subset \bigcap_{j=1}^{\Delta(N)} \Bigl\{\hat z_{jN} \in [l_N,u_N]\Bigr\}\;\;\;\mbox{and}\;\;\;
\lim_{N \to \infty} \Pi_{\theta_N}^N\Bigl(\mathbb{E}_{0N} \cap \bigcap_{j=1}^{\Delta(N)} \Bigl\{\hat z_{jN} \in [l_N,u_N]\Bigr\} \Bigr)=1.
\eqno(F.25)
$$
Indeed, by (ii), $\Delta(N)=o\bigl({\sqrt {N}}/{\log N}\bigr)$. Hence, for $N$ large enough, we have $1/[\Delta(N) +1] > [3 \log N]/\sqrt N$. Therefore, (F.10) and (i) 
imply that 
$$
\log N /\sqrt N \leq J_1(\hat z_{1N}) \leq  J_1(\hat z_{\Delta(N) N}) \leq 1-\log N /\sqrt N. 
$$
By (F.25) we infer 
$$
\Pi_{\theta_N}^N \bigl(\mathbb{E}_{0N} \cap \mathbb{E}_{3N}^c\bigr) \leq 
\Pi_{\theta_N}^N \Bigl(\sup_{z \in [l_N,u_N]} \frac{|\zeta_N(\hat z_{jN})|}{\sqrt{\hat J_N(z)[1-\hat J_N(z)]}} \geq \rho_N[1+o(1)] \Bigr),
$$
where $\rho_N =[\rho_N^{(1)} - \rho_N^{(2)}]/2$ and $\rho_N \asymp \theta_N \sqrt N$. An application of (B.8) finishes the proof. $\hfill{\Box}$\\


\noindent
{\bf Appendix G: Proof of Theorem 6}\\

\noindent
We shall argue that Theorems 2 - 5 imply, via Lemmas 1 and 2, that the regularity assumptions (I.1), (I.2), (II.1) and (II.2) hold true. Besides, (II.1) holds with such $\{\gamma_N\}$ and $\{\lambda_N\}$ that (2.11) is satisfied and Theorem 1 works. 

For ${\cal U}_N^{(I)}={\cal V}_N$ the situation is easy. Theorem 2 implies (I.1) while Theorem 4 along with Lemma 2 yield (I.2).

To verify (II.1) for ${\cal U}_N^{(II)}={\cal T}_N$ it is enough to indicate sequences $\{\gamma_N\}$ and $\{\lambda_N\}$ such that Theorem 3 yields (II.1). Observe that $\gamma_N = \log N$ and $\lambda_N = N/\Delta(N)^{2/(1-\nu)}$ are adequate. Indeed, $Nw_N^2/\lambda_N=[w_N^{1-\nu}\Delta(N)]^{2/{(1-\nu)}} = o(1)$  and Theorem 3 applies with $v=\nu$. Similarly, $Nw_N^2/\gamma_N \to \infty$. Hence (II.1) is proved. 

Assumptions of Theorem 6 are stronger than that of Theorem 5. Therefore, by Theorem 5, (i) and (ii) of Lemma 1 hold true with $ U_2^{(II)}(w)=\Phi(w),\;
U_1^{(II)}(w)=T_1(w)$. Pertaining sequence $\{b_N^{(II)}\}=\{b_{\cal T}(\Pi_{\theta_N}^N)\}$ is of the order $\theta_N \sqrt N$. The distribution function $T_1$ has the property $K(\sqrt{\delta[1-\delta]}w) \leq T_1(w) \leq 2\Phi(w) -1$, where $K(w)=Pr(\sup_{0<t<1}|B(t)| \leq w)$. This implies that $T_1(0)=0$ and, by 
\cite{r102}, $T_1(w)$ is absolutely continuous on $[0,\infty)$. This allows for choosing $w_0^{(II)}$ arbitrarily close to 0 and proves (II.2).

Since the distribution of ${\cal T}_N$ is discrete one and its atoms depend on $N$, therefore (iii) of Lemma 1 deserves some comment. Recall that ${\cal T}_N = \max_{1 \leq j \leq \Delta (N)}
\{-{\cal L}_j\}$; cf. (3.6). Due to stochastic monotonicity of ${\cal T}_N$ one can restrict attention to the case $F=G$. Then the distribution of the vector of ranks is uniform. By (3.4) and (3.5), for each $j=1,...,\Delta (N)$ it holds 
$$
{\cal L}_j = \sqrt{\frac{N}{mn}} \frac{1}{\sqrt{\pi_{jN}(1-\pi_{jN})}} 
\Bigl[\sum_{i=1}^m \frac{n}{N} {\bf 1}_{[0,\pi_{jN})} \Bigl(\frac{R_i - 0.5}{N} \Bigr) - 
\sum_{i=m+1}^N \frac{m}{N} {\bf 1}_{[0,\pi_{jN})} \Bigl(\frac{R_i - 0.5}{N} \Bigr) \Bigr].
$$
Note that the value of ${\cal L}_j$ depends only on the number of $(R_i-0.5)/N,\;i=1,...,m,$ falling into $[0,\pi_{jN})$. Hence, if this number increases by 1 then the first sum in ${\cal L}_j$ increases by $n/N$ while the second one decreases by $m/N$. In consequence, the value of ${\cal L}_j$ increases by 
$$
\delta_{jN} = \sqrt{\frac{N}{mn}} \frac{1}{\sqrt{\pi_{jN}(1-\pi_{jN})}} \leq \delta_N = 4 \sqrt{\frac{\Delta(N)}{N}},\;\;j=1,...,\Delta(N).
$$
Most sparse are locations of atoms of ${\cal L}_1$ and ${\cal L}_{\Delta(N)}$. The minimal value of ${\cal L}_1$ is attained when in the interval $[0,\pi_{1N}]$ ranks of the observations from the first sample are absent. This minimal value, say $L_1$, satisfies
$$-L_1=\lfloor N\pi_{1N}+0.5 \rfloor\frac{m}{N}\sqrt{\frac{N}{mn}}\frac{1}{\sqrt{\pi_{1N}(1-\pi_{1N})}}\asymp
\sqrt{\frac{N\pi_{1N}}{1-\pi_{1N}}}\geq \sqrt{\frac{N}{\Delta(N)}}.$$
Since $\;b_N^{(II)}\asymp \theta_N\sqrt{N}$ the assumption (ii)' yields $\;|L_1/b_N^{(II)}|\to\infty$.

Similar argument applies to ${\cal T}_N$ and yields that the atoms of the distribution of this statistic are located at points with distance not exceeding the above defined $\delta_N$. Hence, in any interval of a fixed length, lying right to the point $b_N^{(II)}=b_{\cal T}(\Pi_{\theta_N}^N)$, there is at least one value of ${\cal T}_N$ and (iii) of Lemma 1 holds.   

Finally, since $b_N^{(II)} \asymp \theta_N \sqrt N$, the assumption (ii)' implies that $[b_N^{(II)}]^2/\lambda_N \to 0$ and $[b_N^{(II)}]^2/\gamma_N \to \infty$ as $N \to \infty$. Therefore, by Lemma 1, (2.11) holds true with $\alpha_N$ given in (2.13). Since $\max_{1 \leq j \leq \Delta(N)}\{\pi_{jN}-\pi_{j-1 N}\}\to 0, $
(2.12) holds, as well, and proves (3.18). $\hfill{\Box}$\\

\vskip 14pt
\noindent {\it Acknowledgements.}
The paper was partially written when B. \'Cmiel was on leave from AGH University of Science and Technology and was granted by postdoc position at the Institute of Mathematics of the Polish Academy of Sciences. 
Moreover, the work of B. \'Cmiel was partially supported by the Faculty of Applied Mathematics AGH UST dean grant for PhD students and young researchers within subsidy of Ministry of Science and Higher Education.
\par

\markboth{\hfill{\footnotesize\rm Tadeusz Inglot, Teresa Ledwina and Bogdan \'Cmiel} \hfill}
{\hfill {\footnotesize\rm INTERMEDIATE EFFICIENCY} \hfill}

\bibhang=1.0pc
\bibsep=1.0pt
\fontsize{9}{14pt plus.8pt minus .6pt}\selectfont

\vskip 14pt
\noindent
Tadeusz Inglot\\
Faculty of Pure and Applied Mathematics, Wroc{\l}aw University of Science and Technology,\\
Wybrze\.ze Wyspia\'nskiego 27, 50-370 Wroc{\l}aw, Poland.
\vskip 2pt
\noindent
E-mail: Tadeusz.Inglot@pwr.edu.pl
\vskip 14pt

\noindent
Teresa Ledwina\\
Institute of Mathematics, Polish Academy of Sciences,\\
ul. Kopernika 18, 51-617 Wroc{\l}aw, Poland.
\vskip 2pt
\noindent
E-mail: ledwina@impan.pl
\vskip 14pt

\noindent
Bogdan \'Cmiel\\
Faculty of Applied Mathematics, AGH University of Science and Technology,\\ 
Al. Mickiewicza 30, 30-059 Cracov, Poland.
\vskip 2pt
\noindent
E-mail: cmielbog@gmail.com
\end{document}